\documentclass[12pt,twoside]{article} 
\usepackage{amssymb,mathrsfs,graphicx}

\newtheorem{theorem}{Theorem}[section]
\newtheorem{lemma}[theorem]{Lemma}
\newtheorem{proposition}[theorem]{Proposition}

\newtheorem{definition}[theorem]{Definition}

\newtheorem{rmrk}[theorem]{Remark}

\makeatletter
\@addtoreset{equation}{section}
\makeatother

\newenvironment{remark}
{\begin{rmrk} \em}
{\end{rmrk}}

\newcommand{\fn} {function}
\newcommand{\me} {measure}
\newcommand{\tr} {trajector}
\newcommand{\erg} {ergodic}
\newcommand{\bi} {billiard}
\newcommand{\sy} {system}
\newcommand{\hyp} {hyperbolic}
\newcommand{\pr} {probability}

\newcommand{\dsy} {dynamical system}

\newcommand{\R} {\mathbb{R}}
\newcommand{\C} {\mathbb{C}}

\newcommand{\Z} {\mathbb{Z}}
\newcommand{\N} {\mathbb{N}}
\newcommand{\qed} {\hfill {\small Q.E.D.} \par\medskip}
\newcommand{\skippar} {\par\medskip}
\newcommand{\ds} {\displaystyle}
\newcommand{\proof} {\noindent \textsc{Proof.} }
\newcommand{\proofof}[1] {\noindent \textsc{Proof of {#1}.} }
\newcommand{\article}[3] {\textsc{{#1}}, {\itshape {#2}}, {{#3}}.}
\newcommand{\book}[3] {\textsc{{#1}}, {\itshape {#2}}, {{#3}}.}
\newcommand{\vol} {\textbf}
\newcommand{\eps} {\varepsilon}
\newcommand{\rset}[2] {\left\{ #1 \: \left| \: #2 \right. \! \right\} }
\newcommand{\lset}[2] {\left\{ \left. \! #1 \: \right| \: #2 \right\} }

\newcommand{\symmdiff} {\triangle}

\renewcommand{\iff} {if and only if\ }

\newcommand{\rw} {random walk}
\newcommand{\ob} {observable}

\newcommand{\ps} {\mathcal{M}}            
\newcommand{\sca} {\mathscr{A}}           
\newcommand{\scb} {\mathscr{B}}           
\newcommand{\scv} {\mathscr{V}}
\newcommand{\G} {\mathbb{G}}              
\newcommand{\ivlim} {\lim_{V \nearrow \ps}}
\newcommand{\go} {\mathcal{G}}            
\newcommand{\lo} {\mathcal{L}}            
\newcommand{\avg} {\overline{\mu}}
\newcommand{\D} {\mathbb{D}}              
\renewcommand{\a} {\alpha}                
\renewcommand{\b} {\beta}                 
\newcommand{\g} {\gamma}                  
\newcommand{\ei}[1] {e^{\imath #1}}       
\newcommand{\emi}[1] {e^{-\imath #1}}     
\newcommand{\wtilde} {\widetilde}
\newcommand{\what} {\widehat}
\newcommand{\wbar} {\overline}
\newcommand{\T} {\mathbb{T}}              
\newcommand{\dpt} {\underline{d} \theta} 
\renewcommand{\t} {\theta}                
\newcommand{\con} {C}                     
\newcommand{\jv} {\mathbb{J}_V}           
\newcommand{\ac} {\mathcal{A}}            

\newenvironment{mylist}
{\begin{list}{}{
\setlength{\labelwidth}{20pt} 
\setlength{\labelsep}{30pt}
\setlength{\leftmargin}{50pt}
\setlength{\rightmargin}{\leftmargin}
\setlength{\itemindent}{0pt}
}}
{\end{list}}

\begin{document}

\title{\textbf{On infinite-volume mixing}}

\author{\textsc{Marco Lenci}
\thanks{
Dipartimento di Matematica, 
Universit\`a di Bologna, 
Piazza di Porta S.\ Donato 5, 
40126 Bologna, Italy.
E-mail: \texttt{lenci@dm.unibo.it} } 
}

\date{Final version for \\
Communications in Mathematical Physics 
\vspace{10pt} \\ 
March 2010}

\maketitle

\begin{abstract}
  In the context of the long-standing issue of mixing in infinite
  ergodic theory, we introduce the idea of mixing for observables
  possessing an infinite-volume average. The idea is borrowed from
  statistical mechanics and appears to be relevant, at least for
  extended systems with a direct physical interpretation.  We discuss
  the pros and cons of a few mathematical definitions that can be
  devised, testing them on a prototypical class of infinite
  measure-preserving dynamical systems, namely, the random walks.

\end{abstract}

\section{Historical introduction}
\label{sec-intro}

The textbook definition of mixing for a transformation $T: \ps
\longrightarrow \ps$ preserving a probability \me\ $\mu$ is 
\begin{equation}
  \label{def-fmix-intro}
  \lim_{n \to \infty} \mu(T^{-n} A \cap B) = \mu(A) \mu(B) 
\end{equation}
for all measurable sets $A, B \subset \ps$ \cite{w}. Extending this
definition to the case where $\mu$ is a $\sigma$-finite \me\ with
$\mu(\ps) = \infty$ is a fundamental issue in infinite \erg\ theory.
References to this problem can be found in the literature at least as
far back as 1937, when Hopf devoted a section of his famous
\emph{Ergodentheorie} \cite{h} to an example of a \dsy\ that he calls
`mixing'. It consists of a set $\ps \subset \R^2$, of infinite Lebesgue
\me, and a map $T: \ps \longrightarrow \ps$ preserving $\mu$, the
Lebesgue \me\ on $\ps$. He proved a property that is equivalent to
this one: there exists a sequence $\{ \rho_n \}_{n \in \N}$ of
positive numbers such that
\begin{equation}
  \label{def-hopf}
  \lim_{n \to \infty} \rho_n \, \mu(T^{-n} A \cap B) = \mu(A) \mu(B)
\end{equation}
for all squarable sets $A,B \subset \ps$ (a squarable set is a bounded
set whose boundary has \me\ zero).

For a long time, the community did not seem to act on this suggestion,
perhaps due in part to the impossibility, in \emph{any} reasonable
\dsy, of verifying (\ref{def-hopf}) for all finite-\me\ sets $A,B$.
(This fact, which might not have been clear to Hopf himself, can be
derived from a famous 1964 paper by Hajian and Kakutani \cite{hk}.)

Work in this direction, however, picked up rather intensively in the
1960's \cite{or, kp, ko, pr, kr, pa}, to the point that Krickeberg in
1967 \cite{kr} proposed (\ref{def-hopf}) as the definition of mixing
for almost-everywhere continuous endomorphisms of a Borel space $(\ps,
\mu)$ (with some extra, inessential, conditions on the sets $A,
B$). Krickeberg applied his definition to Markov chains with an
infinite state space and an infinite invariant \me, which is very
interesting in the context of this paper because our main examples
will be the prototypical infinite-state Markov chains, namely the \rw
s (cf.\ Sections \ref{subs-ex}, \ref{sec-rw} and \ref{sec-pf-main}).

Krickeberg's definition has been studied by several researchers since
then \cite{fr, to1, to2} and, in recent times, it was independently
rediscovered by Isola, who uses (\ref{def-hopf}) with $A = B$ and
calls $\{ \rho_n \}$ the \emph{scaling rate} \cite{i1, i2}. It failed,
however, to establish itself as the ultimate definition of mixing in
infinite \me.  In my opinion, this is not so much because of the
less-than-perfect requirement of a topological structure in a
\me-theoretic problem, but rather for its inherent inability to
describe the ``global'' infinite-\me\ aspects of a dynamics: after
all, (\ref{def-hopf}) only involves finite-\me\ sets. Related to this,
it is unclear how this definition may be specified towards stronger
and more physically relevant chaotic properties, such as, for example,
the rate of correlation decay. A little thinking convinces one that
the speed of convergence in (\ref{def-hopf}) cannot in general be
uniform, even for uniformly nice sets ($A$ and $B$ can be arbitrarily
far from each other so that the l.h.s.\ of (\ref{def-hopf}) is
negligible for arbitrarily long times).

At any rate, by the end of the 1960's, Krengel and Sucheston \cite{ks}
approached the problem from a more \me-theoretic point of view and
devised the following two definitions: A discrete-time, nonsingular
\dsy\ $(\ps, \mu, T)$ is called mixing \iff the sequence $\{ T^{-n} A
\}_{n \in \N}$ is \emph{semiremotely trivial} for all measurable $A
\subset \ps$ with $\mu(A) < \infty$; it is called completely mixing if
the condition holds for \emph{all} measurable $A$. (A nonsingular map
is one for which $\mu(A) = 0$ implies $\mu(T^{-1} A) = 0$. As for the
definition of semiremotely trivial, which is unimportant here, we
refer the reader to \cite{ks}.) Both definitions reduce to the
standard definition (\ref{def-fmix-intro}) for maps preserving a \pr\
\me\ \cite{su}.

However, Krengel and Sucheston themselves proved results that imply
that many reasonable (including all invertible) \me-preserving maps
cannot be completely mixing \cite[Thms.~3.1 and 5.1]{ks}. As for
mixing, again specializing to \me-preserving maps, their definition is
equivalent to
\begin{equation}
  \label{def-tr-mix}
  \lim_{n \to \infty} \mu(T^{-n} A \cap B) = 0
\end{equation}
for all finite-\me\ sets $A, B$ \cite[\S2]{ks}.  This is a rather
brutal weakening of (\ref{def-hopf}): for instance, it would classify
a translation in $\R^d$ as mixing! Therefore, however illuminating
\cite{sa}, the Krengel--Sucheston definitions are of little
applicability to most simple extensive \sy s that mathematicians would
like to study.

Aaronson in 1997 \cite[\S2.5]{a} wrote
\begin{quote}
  [...] the discussion in \cite{ks} indicates that there is no
  reasonable generalisation of mixing. 
\end{quote}
Be that as it may, the drive to produce a general definition of mixing
in infinite \erg\ theory had apparently ceased by the mid 1970's.

\skippar

In this paper I will not try to give a universal and firm definition
of mixing---in which I am not sure I believe myself---for
$\sigma$-finite \me-preserving \dsy s, but rather a few very general
\emph{notions}, that can be completely specified on a case-by-case
basis, depending on what type of information one wants to extract from
the \dsy\ under scrutiny.

To do so, I will borrow some ideas and a little terminology from
physics, in particular from statistical mechanics. The key concept
this work is based on is that of \emph{infinite-volume average}, which
I illustrate with an example.

Let us consider a measurable unbounded $A \subset \R^2$ and ask, what
is the probability that a random $x \in \R^2$ belongs to $A$? Clearly,
the answer fully depends on what we mean by random.  Suppose we
specify that random means that each $x$ can be drawn with equal
probability.  Then the question itself no longer makes sense because
the Lebesgue \me\ $m$ on $\R^2$, which is the only uniform \me\ on
$\R^2$, cannot be normalized.

However, remembering that long-gone course in statistical mechanics,
one might come up with the idea that the sought probability is
something like
\begin{equation}
  \label{avg-ex}
  \lim_{r \to +\infty} \frac{ m( A \cap [-r,r]^2 ) } {4r^2},
\end{equation}
provided the limit exists.  Of course, such an answer is riddled with
issues, but it does capture the idea that in physics one only looks at
finite quantities. Infinity is a mental construct to fit an endless
amount of situations, and a finite limit at infinity is the formal way
to say that most of these situations will look alike.

More generally, for a given \dsy\ $(\ps, \sca, \mu, \{ T^t \} )$,
where $(\ps, \sca)$ is a \me\ space and $\{ T^t \}$ is a (semi-)group
of transformations $\ps \longrightarrow \ps$ preserving the infinite
\me\ $\mu$, we will choose a family of ever-larger sets $V$, with
$\mu(V) < \infty$, that ``approximate $\ps$''.  Using the language of
physics, one might say that choosing these sets will define how we
\me\ our infinite \sy---more precisely, how we \me\ its \emph{\ob s}.

We will deal with two types of \ob s: The \emph{global}, or
\emph{macroscopic}, \ob s will be a suitable class of \fn s $F \in
L^\infty (\ps, \sca, \mu)$ for which
\begin{equation}
  \label{avg-intro}
  \avg{(F)} := \ivlim \, \frac1{\mu(V)} \int_V F \, d\mu
\end{equation}
exists. What the above limit means and what class of \fn s it applies
to will be clarified in Section \ref{sec-iv}.  The \emph{local}, or
\emph{microscopic}, \ob s will be essentially the elements of $L^1
(\ps, \sca, \mu)$.

Then, skipping many necessary details and all-important specifications
which are found in Sections \ref{sec-iv} and \ref{sec-defs}, our
notions of mixing will basically reduce to the two limits:
\begin{equation}
  \label{ggmix}
  \lim_{t \to \infty} \avg((F \circ T^t) G) = \avg(F) \avg(G),  
\end{equation}
for any two global \ob s $F, G$; and
\begin{equation}
  \label{glmix}
  \lim_{t \to \infty} \mu((F \circ T^t) g) = \avg(F) \mu(g),
\end{equation}
for $F$ global and $g$ local (with the obvious notation $\mu(g) :=
\int g \, d\mu$). 

To the extent to which the above notions can be made into rigorous
definitions---and they can, cf.\ Section \ref{sec-defs}---they seem to
improve on the attempted definitions that we have recalled earlier,
chiefly because they involve \ob s which can be supported throughout
the phase space (think, for example, of the velocity of a particle in
an aperiodic Lorentz gas, or the potential energy of a small mass in a
formally infinite celestial conglomerate, etc.). So they are more apt
to reveal the large-scale aspects of a given dynamics.

In particular, (\ref{ggmix}) may be called \emph{global-global
mixing}, because it somehow expresses the vanishing of the
correlation coefficient between two global \ob s, while (\ref{glmix})
may be called \emph{global-local mixing}, because the coupling is
between a global and a local \ob. The latter notion can be quite
useful if one takes $g \ge 0$ with $\mu(g) = 1$. Then $\mu_g$, the
probability \me\ defined by $d\mu_g := g \, d\mu$, can be considered
an initial state for the \sy. In this interpretation, the l.h.s.\ of
(\ref{glmix}) reads $\mu_g (F \circ T^t) =: T_*^t \mu_g (F)$, that is,
the expected value of the \ob\ $F$ relative to the state at time $t$,
and (\ref{glmix}) asserts that such quantity converges to
$\avg(F)$. Hence, $\avg$ acts as a sort of equilibrium state for the
\sy.

In Section \ref{sec-rw} we will apply the above ideas to certain basic
yet representative examples of infinite-\me\ \dsy s, the \rw s in
$\Z^d$. We will specify all the mathematical objects needed to obtain
rigorous definitions out of (\ref{ggmix})-(\ref{glmix}) and we will
check that, under reasonable conditions, all these definitions are
verified. The proofs, given in Section \ref{sec-pf-main}, use basic
harmonic analysis on groups \cite{r} and an estimate for a certain
Fourier norm. The latter is presented in the Appendix, together with
other technical results.

\bigskip

\noindent
\textbf{Acknowledgments.}\ I would like to thank an anonymous referee
for pointing out a relevant mistake in the first draft of the
manuscript. 

\section{Infinite-volume limit and observables}
\label{sec-iv}

A \emph{\me-preserving \dsy} is the quadruple $(\ps, \sca, \mu,
\{T^t\} )$, where $\ps$ is a \me\ space with the $\sigma$-algebra
$\sca$, endowed with the $\sigma$-finite \me\ $\mu$, while $\{ T^t
\}_{t \in \G}$ is a group (respectively, semigroup) of automorphisms
(respectively, endomorphisms) $\ps \longrightarrow \ps$, labeled by
the free additive parameter $t \in \G$ (without significant loss of
generality, we assume $\G = \N$, $\Z$, or $\R$).  This means that
\begin{equation}
  \label{meas-pres} 
  \mu(T^{-t} A) = \mu(A), \qquad \forall A \in \sca, \forall t \in \G.
\end{equation}
If $\G = \R$, $\{ T^t \}$ is called \emph{the flow}, whereas if $\G =
\Z$ or $\N$, the generator $T := T^1$ is called \emph{the map}, and
one usually denotes the \dsy\ by $(\ps, \mu, T)$.

In this paper we are interested in the case $\mu(\ps) = \infty$. The
\me-theoretic properties of \dsy s preserving an infinite \me\ are the
subject of \emph{infinite \erg\ theory} \cite{a}, whose most basic
definition, perhaps, is the following \cite{hk}:
\begin{definition}
  \label{def-erg}
  The \me-preserving \dsy\ $(\ps, \sca, \mu, \{ T^t \} )$ is called
  \emph{\erg} if every measurable invariant set (i.e., any $A \in
  \sca$ such that $\mu( T^{-t} A \,\symmdiff\, A) = 0$ $\forall t$), has
  either zero \me\ or full \me\ (the latter meaning $\mu (\ps
  \setminus A) = 0$).
\end{definition}
In the case where $\mu$ is a \pr\ \me, the above is one of the several
equivalent formulations of \erg ity, including among others: the
equivalence of the Birkhoff average with the phase average, for all $f
\in L^1(\ps, \sca, \mu)$ (the definition generally ascribed to
Boltzmann); the absence of nontrivial integrals of motion in $L^1(\ps,
\sca, \mu)$; the strong law of large numbers for the random variables
$\{ f \circ T^t \}$. In infinite \me, these definitions are no longer
equivalent and, among those that keep making sense, Definition
\ref{def-erg} is in some sense the strongest.

A notion that is much harder to transport to infinite \erg\ theory, as
we have discussed in the introduction, is that of \emph{mixing}.  In
terms of \emph{\ob s}, i.e., scalar \fn s on $\ps$, it reads:
\begin{equation}
  \label{def-fmix}
  \forall f,g \in L^2(\ps, \sca, \mu), \qquad \lim_{t \to \infty} 
  \mu((f \circ T^t) g) = \mu(f) \mu(g). 
\end{equation}
In other words, the correlation coefficient between the two random
variables $f \circ T^t$ and $g$ vanishes asymptotically.  This
phrasing makes it apparent that the notion of finite-\me\ mixing is
intrinsically probabilistic.

\subsection{Infinite-volume limit}
\label{subs-iv}

The discussion in the introduction suggests that one should find an
asymptotic decorrelation formula, similar to (\ref{def-fmix}), which
applies to \ob s that, unlike $L^2$ \fn s, ``see'' a nonnegligible
portion of the space.  For this, one needs to define a sort of
``normalized \me'' for these \ob s.  This is how one comes to think of
averaging a \fn\ over $\ps$ by means of an \emph{infinite-volume
limit}.

The idea is borrowed from statistical mechanics where the question
arises of measuring the \emph{intensive} quantities of a \sy\ (for
example, the temperature or the density of a gas).  These are
represented by sequences of \fn s defined over larger and larger phase
spaces, corresponding to larger and larger portions of the physical
\sy, normalized in such a way that their integral converges to a
finite number. This number is supposed to predict the outcome of an
experimental measurement.

\skippar

Coming back to our math, we introduce a notation that is going to be
used often in the remainder:
\begin{equation}
  \label{af}
  \sca_f := \rset{A \in \sca} {\mu(A) < \infty}.
\end{equation}

\begin{definition}
  \label{def-exhau}
  The family $\scv \subset \sca_f$ is called \emph{exhaustive} if it
  contains a sequence $\{ V_n \}_{n \in \N}$, increasing w.r.t.\
  inclusion, such that $\bigcup_n V_n = \ps$.
\end{definition}

\begin{definition}
  \label{def-iv}
  Let $\scv$ be exhaustive.  If $\phi$ is defined on $\sca_f$ and has
  values in some topological space, we write
  \begin{displaymath}
    \ivlim \phi(V) \: := \lim_{{V \in \scv} \atop {\mu(V) \to \infty}}
    \phi(V) = L,
  \end{displaymath}
  when, for every neighborhood $\mathcal{U}$ of $L$, there exists an
  $M \in \R^+$ such that 
  \begin{displaymath}
    V \in \scv, \mu(V) \ge M \ \Longrightarrow\ \phi(V) \in \mathcal{U}.
  \end{displaymath}
  This will be called the \emph{$\mu$-uniform infinite-volume limit
  w.r.t.\ the family $\scv$}, or simply the \emph{infinite-volume
  limit}.
\end{definition}

It is apparent that the above definition depends decisively on the
choice of $\scv$. In the example discussed in the introduction, cf.\
(\ref{avg-ex}), $\scv = \{ [-r,r]^2 \}_{r>0}$, and it is easy to think
of other exhaustive families of subsets of $\R^2$ for which the
infinite-volume limit for the ``\pr'' of $A$ differs from
(\ref{avg-ex}).

In general, all the results that we are going to discuss in this paper
will depend on the choice of $\scv$ in Definition \ref{def-iv}. This
is no shortcoming! In fact, we \emph{want} to retain this choice,
because this is how we incorporate in the mathematical description of
an extended \sy\ the way we \emph{observe} the \sy, that is, how we
\me\ the \ob s of that \sy. In other words, the choice of $\scv$
defines what it means to pick a large region of $\ps$, which we
assume---or rather \emph{declare}---represents the whole space.

The first property we assume of our \sy s may be called `compatibility
of the infinite-volume limit with the dynamics':
\begin{mylist}
\item[\textbf{(A1)}] For any fixed $t \in \G$, for $V \nearrow \ps$, \
  $\ds \mu( T^{-t} V \symmdiff V ) = o(\mu(V))$.
\end{mylist}
This means that the scale of the dynamics is local and not global: If
$V$ is a very large set ``approximating $\ps$'', then its evolution at
a fixed time should differ little from $V$, in relative terms. (In
statistical mechanics one would say that, over a large region of the
space, the dynamics can only produce negligible \emph{surface
effects}---we will come back to this point in Section
\ref{subs-ggm}.)

One can expect \textbf{(A1)} to hold in most situations. It does, for
instance, when $\ps$ is a metric space such that the $\mu$-measure of a
ball grows like a power of the radius, uniformly in the center of the
ball; the dynamics is bounded, i.e., $\forall t \in \G$, $\exists
K=K(t)$ such that, for $\mu$-a.e.\ $x \in \ps$, $\mathrm{dist} (T^t x,
x) \le K$; and the elements of $\scv$ are balls.

\subsection{A couple of examples}
\label{subs-ex}

A very simple example of this setup is the \dsy\ defined as
follows. Set
\begin{equation}
  \label{exrw1}
  \varphi(x) := \left\{ 
  \begin{array}{lll}
    3x - 1, && \mbox{ for } x\in [0,1); \\
    0 , && \mbox{ otherwise.}
  \end{array}
  \right.
\end{equation}
Then, on $\ps = \R$, define the self-map
\begin{equation}
  \label{exrw2}
  Tx := \sum_{j \in \Z} \varphi(x-j) + j.
\end{equation}
The definition is well-posed because, for any given $x \in \R$, only
one term of the sum is nonzero. Furthermore, each $x$ has three
distinct counterimages via $T$, where the derivative of the map is
constantly equal to 3. This shows that $T$ is a noninvertible map
which preserves the Lebesgue \me\ $\mu$.

$(\ps, \mu, T)$ describes a \rw\ in $\Z$, in the following sense:
Suppose an initial condition $x \in [0,1)$ is randomly chosen
according to $\mu_0$, the Lebesgue \me\ in $[0,1)$ (it is no loss of
generality to restrict to $[0,1)$ because the \dsy\ is clearly
invariant for the action of $\Z$). Then $T x$ will land in one of the
intervals $[-1,0)$, $[0,1)$, $[1,2)$, with \pr\ 1/3 in each case. More
generally, if $[ x ] := \max \rset{m \in \Z} {m \le x}$ denotes the
integer part of $x \in \R$, we have
\begin{equation}
  \mu_0 \left( \rset{x \in [0,1) } { \, [T^n x] = k_n, [T^{n-1} x] =
  k_{n-1}, \dots, [Tx] = k_1 } \right) = 3^{-n},
\end{equation}
provided that $|k_j - k_{j-1}| \le 1$, for $j=1, \ldots, n$ (with
$k_0=0$). This implies that, using the notation of conditional \pr,
\begin{eqnarray}
  && \mu_0 \left( \, [T^n x] = k_n \: \left| \ [T^{n-1} x] = k_{n-1},
  \dots, [Tx] = k_1 \right. \right) = \nonumber \\
  && \mu_0 \left( \, [T^n x] = k_n \: \left| \ [T^{n-1} x] = k_{n-1} 
  \right. \right) = 1/3.
\end{eqnarray}
Hence, $\{ [T^n x] \}_n$ is a Markov chain in $\Z$ with same-site
and nearest-neighbor jumps, each with \pr\ 1/3; namely, it is a
(space-)homogeneous \rw.

As for the choice of $\scv$, the example of the introduction would
seem to suggest that we pick sets of the type $V = [-r,r]$ (with $r
\in \R$) or, to fully exploit the $\Z$-structure of this \dsy\ with no
appreciable loss of generality, sets of the type $V = [-k,k]$ (with $k
\in \N$). Although this is a legitimate choice, we will see later that
a better option is
\begin{equation}
  \label{good-scv-rw}
  \scv := \rset{ [k, \ell] } {k,\ell \in \Z, k<\ell}.
\end{equation}
(Actually, this will be a crucial part of our discussion, and we refer
the reader to Section \ref{subs-ggm}.)  For $V = [k,\ell]$, $n \in
\N$, and $\ell-k > 2n$, it is easy to verify that
\begin{equation}
  \label{pf-a1-sec2}
  [k+n, \ell-n] \subset T^{-n} V \subset [k-n, \ell+n],
\end{equation}
which implies \textbf{(A1)}.

\skippar

A less trivial \sy\ that fits well the framework we are describing is
the (aperiodic) Lorentz gas \cite{l1, l2}.  In $\R^2$ (just to fix the
smallest interesting dimension) a Lorentz gas is the \bi\ \sy\ in $C
:= \R^2 \setminus \bigcup_{n \in \N} \mathcal{O}_n$, where $\{
\mathcal{O}_n \}$ is a countable collection of pairwise disjoint,
convex, bounded regular sets. This means that a point particle moves
with constant unit velocity in $C$ until it hits an obstacle
$\mathcal{O}_n$, which reflects the particle according to the Fresnel
law: the angle of reflection equals the angle of incidence (the
modulus of the velocity remains equal to 1). The phace space for this
\sy\ is then $\ps = C \times S^1$, where $q \in C$ represents the
position and $v \in S^1$ the velocity of the particle. If $T^t$
denotes the flow just described (which is unambiguously defined at all
noncollision times), it is well known that $T^t$ preserves the
Liouville \me\ $\mu$, which turns out to be the product of the
Lebesgue \me\ on $C$ and the Haar \me\ on $S^1$. Clearly, save for
bizarre situations, $\mu(\ps) = \infty$.

Since sufficient conditions for the \erg ity of $(\ps, \mu, T)$ are
known \cite{l1}, and since, for the finite-\me\ version of the Lorentz
gas (the so-called \emph{Sinai \bi}), mixing and stronger stochastic
properties essentially follow from \erg ity \cite{si, bs, cm}, it is
of interest to devise one or more sound definitions of mixing for this
\dsy.

As for the exhaustive family $\scv$, in analogy to the previous \sy, a
reasonable choice would be
\begin{equation}
  \label{good-scv-lg}
  \scv := \rset{ (C \cap R) \times S^1 } {R = [a,b] \times [c,d],
  \mbox{ with } a<b,\, c<d}.
\end{equation}

In Section \ref{sec-rw} we will present a third example, which
generalizes, in more than one way, the first \sy\ introduced above. It
is a class of invertible \dsy s representing all homogeneous \rw s in
$\Z^d$. One of its points of relevance is that it is designed to
retain the most essential features of the Lorentz gas discussed
above. It is thus a greatly simplified toy model, which we are able to
study in depth. As a matter of fact, it will be the testing ground for
our new notions of mixing, cf.\ Sections \ref{sec-rw} and
\ref{sec-pf-main}.

\subsection{Global and local observables}
\label{subs-obs}

In order to define a surrogate \pr\ \me\ for our \sy, we need to
declare \emph{what} we intend to \me. In other words, we need to
specify the \emph{\ob s}, namely, the \fn s $\ps \longrightarrow \R$
which represent the (sole) information that we can get on the state of
the \sy.

In finite \erg\ theory, this class of \fn s is $L^1(\ps, \sca, \mu)$,
or sometimes $L^2(\ps, \sca, \mu)$. In virtually every situation, both
are amply sufficient to give a full description of the state of the
\sy\ (in fact, quite generally, the position itself of $x \in \ps$ is
given by a finite number of square-integrable \fn s).  This is
conspicuously not true in infinite \erg\ theory. Indeed, the
forthcoming discussion will try to convince the reader that the choice
of the \ob s is precisely at the heart of the matter in infinite-\me\
mixing (a point that, after all, was already contained in \cite{ks}).

We deal with two categories of \ob s: the \emph{global}, or
\emph{macroscopic}, \ob s and the \emph{local}, or \emph{microscopic},
\ob s.

Let us start by introducing the former, whose space we denote by
$\go$.  We will not give a definition, but rather a presentation of
the minimal features that $\go$ should have. A precise definition only
makes sense on a case-by-case basis and, indeed, the choice of $\go$
is part of our description of the \sy, just like the choice of $\scv$.
As a typographical rule, we indicate a global \ob\ with an upper-case
Roman letter, as in $F: \ps \longrightarrow \R$.

We require at least the following conditions:
\begin{mylist}
\item[\textbf{(A2)}] $\ds \go \subset L^\infty (\ps, \sca, \mu)$.

\item[\textbf{(A3)}] $\ds \forall F \in \go, \qquad \exists \, \avg(F)
  := \ivlim\, \mu_V (F) := \ivlim\, \frac1 {\mu(V)} \int_V F \, d\mu$.
\end{mylist}
We call $\avg(F)$ the \emph{average} of $F$ (w.r.t.\ $\mu$ and
$\scv$). This functional is dynamics-invariant:

\begin{lemma}
  \label{lem-inv}
  Under assumptions {\bf (A1)-(A3)}, $\avg(F) = \avg(F \circ T^t)$,
  $\forall t \in \G$.
\end{lemma}

\proof Using the invariance of $\mu$ and then \textbf{(A1)-(A2)}, we
have
\begin{equation}
  \frac1 {\mu(V)} \int_V F \, d\mu = \frac1 {\mu(V)} \int_{T^{-t} V} 
  (F \circ T^t) \, d\mu = \frac1 {\mu(V)} \int_V (F \circ T^t) \, 
  d\mu + o(1).
\end{equation}
Applying \textbf{(A3)} gives the assertion.
\qed

As for the class of local \ob s, denoted by $\lo$, this can be
generally taken to be $L^1(\ps, \sca, \mu)$. As we will see below,
this choice is much less delicate than the choice of
$\go$. Nonetheless, some results may require occasional restrictions
on $L^1$, so, in the same spirit as \textbf{(A2)}, we only require
\begin{mylist}
\item[\textbf{(A4)}] $\ds \lo \subseteq L^1 (\ps, \sca, \mu)$.
\end{mylist}
Local \ob s are indicated with a lower-case Roman letter, as in $g:
\ps \longrightarrow \R$.

\section{Definitions and related questions}
\label{sec-defs}

\subsection{Global-global mixing}
\label{subs-ggm}

On the basis of Lemma \ref{lem-inv}, one might attempt the following
definition of mixing:
\begin{mylist}
\item[\textbf{(M1)}] $\ds \forall F,G \in \go, \qquad \lim_{t \to
    \infty} \, \avg ( (F \circ T^t) G ) = \avg(F) \, \avg(G)$,
\end{mylist}
provided that $\avg( (F \circ T^t) G )$ exists for all $t \in \G$, or
at least for $t$ large enough. This last point represents a problem,
because it is not easy, in general, to devise a space $\go$ with the
property that $F, G \in \go$ implies $(F \circ T^t) G \in \go$ for all
large $t$ (sometimes $\go$ is not even $T^t$-invariant, cf.\
(\ref{def-gom}) later on).

Generally speaking, there are only two solutions to this
problem---which would be more honestly described as ways around
it. The first solution is to declare that this question should be
dealt with on a case-by-case basis.  The second solution is to devise
another definition of mixing which just does away with the problem:
\begin{mylist}
\item[\textbf{(M2)}] $\ds \forall F,G \in \go, \qquad \lim_{{t \to
  \infty} \atop {V \nearrow \ps}} \, \mu_V ( (F \circ T^t) G ) =
  \avg(F) \, \avg(G)$,
\end{mylist}
having adopted the notation $\mu_V (\cdot) = \int_V (\cdot) d\mu /
\mu(V)$, as introduced in \textbf{(A3)}. The above means that,
$\forall \eps > 0$, $\exists M > 0$ such that, for all $t \ge M$ and
$V \in \scv$ with $\mu(V) \ge M$, 
\begin{equation}
  \label{defs-05}
  \left| \mu_V ( (F \circ T^t) G ) - \avg(F) \avg(G) \right| < \eps.
\end{equation}

Though cast in a less polished form than \textbf{(M1)}, \textbf{(M2)}
still retains a great deal of the physical meaning of mixing because
it prescribes that, if the region $V$ is big enough and the time $t$
is large enough, the two \ob s $F \circ T^t$ and $G$ are practically
uncorrelated on $V$. Actually, in some sense, \textbf{(M2)} is even
stronger than \textbf{(M1)}, because it implies that, fixed $V$,
(\ref{defs-05}) occurs \emph{uniformly} in $t$, for $t$ large. The
same is not guaranteed by \textbf{(M1)}. See also Proposition
\ref{prop-hye} later on.

We refer to \textbf{(M1)} and \textbf{(M2)} as definitions of
\emph{global-global mixing} because they consider the coupling of two
global \ob s.

\skippar

Let us now focus on a couple of less technical and more substantial
questions concerning both \textbf{(M1)} and \textbf{(M2)}. The first
has to do with the importance on $\scv$ too, not just $\go$, for
either condition to \fn\ as a sound definition of mixing.

Let us exemplify the question by means of the \dsy\ defined by
(\ref{exrw1})-(\ref{exrw2}). This is a \sy\ that should be classified
as mixing by any reasonable definition (cf.\ also Section
\ref{sec-rw}). Suppose that for that \sy\ we had made the first, more
restrictive, choice of $\scv$ presented in Section \ref{subs-ex}, that
is, we had chosen sets of the type $V = [-k,k]$, with $k \in \N$. The
\fn\ $F: \R \longrightarrow \R$, defined by
\begin{equation}
  \label{defs-10}
  F(x) := \left\{
  \begin{array}{rll}
    -1, && \mbox{for } x < 0 ; \\
    1, && \mbox{for } x \ge 0 ,
  \end{array} \right.
\end{equation}
is bounded and has average $\avg(F) = 0$.  Now, fix $n > 0$ and
consider $F \circ T^n$. Given the action of the map $T$ and its
interpretation as a \rw, it is not hard to see that, for $x < -n$,
$F(T^n x) = -1$ and, for $x \ge n$, $F(T^n x) = 1$ (determining $F(T^n
x)$ for $x \in [-n, n)$ is more complicated and irrelevant here). This
and (\ref{defs-10}) imply that, for $|x| > n$, $F(T^n x) F(x) =
1$. Therefore, for $k$ much larger than $n$ and $V = [-k,k]$, $\mu_V (
(F \circ T^n) F)$ is very close to 1 and indeed $\avg( (F \circ T^n)
F) = 1$. Since $\avg(F) = 0$, this shows that, for $F$ as in
(\ref{defs-10}) and $G=F$, both limits in \textbf{(M1)} and
\textbf{(M2)} fail to hold!

But this is reasonable: after all, $F$ has variations (causing it to
be nonconstant) only on a negligible set, namely $\{ x=0 \} \subset
\ps$. By negligible set we mean, in this context, a set whose
$\rho$-neighborhoods, for all $\rho>0$, have ``\me'' zero w.r.t.\
$\avg$. Therefore this is an instance of the phenomenon, which is well
known in statistical mechanics, whereby the infinite-volume limit does
not see \emph{surface effects}. (When the dynamics is bounded, in the
sense specified in the last paragraph of Section \ref{subs-iv}, the
evolution of $F$ can produce no more than surface effects.)

So we must avoid global \ob s that have significant variations on
negligible sets. But this does not mean that we should cherry-pick our
\ob s (although there is nothing wrong with that)! In the case at
hand, for instance, the unusable \ob s are automatically eliminated by
a smarter choice of $\scv$, the one given in (\ref{good-scv-rw}). That
exhaustive family is translation invariant, which seems right for a
\sy\ that is translation invariant.

An analogous discussion can be made for the other example presented in
Section \ref{subs-ex}, the Lorentz gas, and for most extended \dsy s
one can imagine. 

We may conclude that $\scv$ should not be so small as to make the
verification of \textbf{(M1)} impossible nor, at the same time, so
large as to make the class of global \ob s satisfying \textbf{(A3)}
too meager. A happy medium might be for $\scv$ to include all the
symmetries, or ``quasi-symmetries'', of the \sy\ and no more.

\skippar

The second, and more critical, question concerning \textbf{(M1)-(M2)}
is that these definitions are completely blind to the local aspects of
the dynamics. For instance, they are not able to detect an invariant
set $A$, if $\mu(A) < \infty$. One can easily produce a \sy\ that is
mixing in one of the above senses, but not \erg\ as per Definition
\ref{def-erg}. (For example, take a Lorentz gas and make one scatterer
hollow: points inside that scatterer will stay confined there, thus
breaking \erg ity, while all other \tr ies will be the same as in the
unperturbed \sy, which one believes to be at least
\textbf{(M1)}-mixing, with the right choice of $\go$ and $\scv$.)

We have no fix for this issue, other then giving a few more
definitions which take into account local \ob s as well.

\subsection{Global-local mixing}
\label{subs-glm}

The most natural way to couple global and local \ob s in a definition
of mixing is this:
\begin{mylist}
\item[\textbf{(M4)}] $\ds \forall F \in \go, \ \forall g \in \lo,
  \qquad \lim_{t \to \infty} \, \mu((F \circ T^t) g) = \avg(F)
  \mu(g)$,
\end{mylist}
where we have used the convenient notation $\mu(g) := \int_\ps g \,
d\mu$.  This is the first notion of \emph{global-local mixing} we
give. Since for some \sy s, such as the Lorentz gas, this can be
rather hard to prove \cite{l3}, we give a weaker version as well:
\begin{mylist}
\item[\textbf{(M3)}] $\ds \forall F \in \go, \ \forall g \in \lo$ with
  $\ds \mu(g)=0, \qquad \lim_{t \to \infty} \, \mu((F \circ T^t) g) =
  0$.
\end{mylist}
(Notice that \textbf{(M3)} and \textbf{(M4)} are equivalent in
ordinary \erg\ theory, because one can always subtract a constant
\fn\ from any \ob\ in order to make its integral vanish. Not so in
infinite \me!)

We will see momentarily that, for \sy s for which a uniform version
of \textbf{(M4)} can be established, it is possible to pass from
global-local mixing to global-local mixing.  So we give one last
definition:
\begin{mylist}
\item[\textbf{(M5)}] $\ds \forall F \in \go, \qquad \lim_{t \to
  \infty} \, \sup_{g \in \lo \setminus 0 } \: \frac1 {\mu(|g|)}
  \left| \mu((F \circ T^t) g) - \avg(F) \mu(g) \right| = 0$.
\end{mylist}

\subsection{Summary of assumptions and definitions}
\label{subs-summ}

For the convenience of the reader, we summarize here all the
assumptions we have made and all the definitions of mixing we have
given, listing the latter in the correct hierarchical order, as
clarified by Propositions \ref{prop-hye} and \ref{prop-m5m2} below.

\skippar

The following are the minimal requirements on the \dsy\ $(\ps, \sca,
\mu, \{ T^t \} )$, the exhaustive family $\scv$, the space of the
global \ob s $\go$, and the space of the local \ob s $\lo$:

\begin{mylist}
\item[\textbf{(A1)}] For any fixed $t \in \G$, for $V \nearrow \ps$, \
  $\ds \mu( T^{-t} V \symmdiff V ) = o(\mu(V))$.

\item[\textbf{(A2)}] $\ds \go \subset L^\infty (\ps, \sca, \mu)$.

\item[\textbf{(A3)}] $\ds \forall F \in \go, \qquad \exists \, \avg(F)
  := \ivlim\, \mu_V (F) := \ivlim\, \frac1 {\mu(V)} \int_V F \, d\mu$.

\item[\textbf{(A4)}] $\ds \lo \subseteq L^1 (\ps, \sca, \mu)$.
\end{mylist}

The definitions of global-global mixing are:

\begin{mylist}
\item[\textbf{(M1)}] $\ds \forall F,G \in \go, \qquad \lim_{t \to
  \infty} \, \avg ( (F \circ T^t) G ) = \avg(F) \, \avg(G)$.

\item[\textbf{(M2)}] $\ds \forall F,G \in \go, \qquad \lim_{{t \to
  \infty} \atop {V \nearrow \ps}} \, \mu_V ( (F \circ T^t) G ) =
  \avg(F) \, \avg(G)$.
\end{mylist}

The definitions of global-local mixing are:

\begin{mylist}
\item[\textbf{(M3)}] $\ds \forall F \in \go, \ \forall g \in \lo$ with
  $\ds \mu(g)=0, \qquad \lim_{t \to \infty} \, \mu((F \circ T^t) g) =
  0$.

\item[\textbf{(M4)}] $\ds \forall F \in \go, \ \forall g \in \lo,
  \qquad \lim_{t \to \infty} \, \mu((F \circ T^t) g) = \avg(F)
  \mu(g)$.

\item[\textbf{(M5)}] $\ds \forall F \in \go, \qquad \lim_{t \to
  \infty} \, \sup_{g \in \lo \setminus 0 } \: \frac1 {\mu(|g|)}
  \left| \mu((F \circ T^t) g) - \avg(F) \mu(g) \right| = 0$.
\end{mylist}

\begin{proposition}
  \label{prop-hye}
  Under all the assumptions made so far,
  \begin{displaymath}
    \mathbf{(M5)} \ \Longrightarrow \ \mathbf{(M4)} \ \Longrightarrow 
    \ \mathbf{(M3)}. 
  \end{displaymath}
  Furthermore, \emph{\textbf{(M2)}} implies that the limit in
  \emph{\textbf{(M1)}} holds for all pairs $F,G \in \go$ such that
  $\avg((F \circ T^t)G)$ exists for all $t$ large enough.
\end{proposition}

\proof The chain of implications is obvious. The last assertion
follows directly from the definition of the double limit $t \to
\infty$, $V \nearrow \ps$; cf.\ Section \ref{subs-ggm}.
\qed

With reasonable hypotheses, the strongest version of global-local
mixing implies the ``strongest'' version of global-global mixing:

\begin{proposition}
  \label{prop-m5m2}
  Suppose that every $G \in \go$ can be written $\mu$-almost
  everywhere as
  \begin{displaymath}
    G(x) = \sum_{j \in \N} g_j(x), \qquad \mbox{with } g_j \in \lo,
  \end{displaymath}
  and, for every $V \in \scv$, there exists a finite subset $\jv$ of
  $\N$, such that
  \begin{eqnarray}
    \label{m5m2-h1}
    \mu\left( \left| \, G \, 1_V - \sum_{j \in \jv} g_j \right|
    \right) &=& o(\mu(V)) ; \\
    \label{m5m2-h2}
    \sum_{j \in \jv} \| g_j \|_{L^1} &=& O(\mu(V)) . 
  \end{eqnarray}
  Then \emph{\textbf{(M5)}} $\Longrightarrow$ \emph{\textbf{(M2)}}.
\end{proposition}

\begin{remark}
  \label{rk-m5m2}
  The hypotheses of Proposition \ref{prop-m5m2} above are less
  cumbersome than they appear. One should think of the very common
  situation in which $\ps$ admits a partition of unity, $\sum_j
  \psi_j(x) \equiv 1$, where the $\psi_j$ are nonnegative integrable
  \fn s which are roughly translations of one another. In many such
  cases one can expect $g_j := G \psi_j$ to verify all of the above
  conditions. At any rate, if the $g_j$ are all nonnegative or all
  nonpositive, then (\ref{m5m2-h2}) follows from (\ref{m5m2-h1}).
\end{remark}

\proofof{Proposition \ref{prop-m5m2}} Fix $F,G \in \go$. We may assume
that $\avg(F) \ne 0$, otherwise in the following argument we consider
$F_c := F+c$, where $c$ is a nonnull constant, and easily derive the
sought result at the end of the proof.

Take $\eps>0$ and denote for short $F^t := F \circ T^t$ and $g_V :=
\sum_{j \in \jv} g_j$. (\ref{m5m2-h1}) and \textbf{(A3)} imply that,
for $\mu(V)$ large enough,
\begin{equation}
  \label{m5m2-10}
  \left| \frac{\mu(g_V)} {\mu(V)} - \avg(G) \right| \le \left| 
  \frac{\mu(g_V)} {\mu(V)} - \mu_V (G) \right| + \left| \mu_V (G) - 
  \avg(G) \right| \le \frac{\eps} { 3 \, |\avg(F)| }
\end{equation}
and, since $F$ is bounded,
\begin{equation}
  \label{m5m2-20}
  \left| \mu_V (F^t G) - \frac{\mu(F^t g_V)} {\mu(V)} \right| \le 
  \frac{\eps}3.
\end{equation}
On the other hand, \textbf{(M5)} implies that
\begin{equation}
  \label{m5m2-30}
  \left| \mu(F^t g_j) - \avg(F) \mu(g_j) \right| \le \vartheta(t) \|
  g_j \|_{L^1} ,
\end{equation}
where $\ds \lim_{t \to \infty} \vartheta(t) = 0$ and $\vartheta(t)$
does not depend on $j$ or $V$. Summing over $j \in \jv$ and using
(\ref{m5m2-h2}), one gets
\begin{equation}
  \label{m5m2-40}
  \left| \frac{\mu(F^t g_V)} {\mu(V)} - \frac{ \avg(F) \mu(g_V) }
  {\mu(V)} \right| \le \vartheta(t) \, \frac{ \sum_{j \in \jv} 
  \| g_j \|} {\mu(V)} \le \frac{\eps}3,
\end{equation}
for both $\mu(V)$ and $t$ large enough.

Putting together (\ref{m5m2-20}), (\ref{m5m2-40}) and (\ref{m5m2-10}),
in that order, we conclude that there exists $M = M(\eps)$ such that,
for $\mu(V) \ge M$ and $t \ge M$,
\begin{equation}
  \label{m5m2-50}
  \left| \mu_V (F^t G) - \avg(F) \avg(G) \right| \le \eps,
\end{equation}
which is precisely \textbf{(M2)}.
\qed 

\section{Mixing for random walks}
\label{sec-rw}

In this section we see how the previous definitions play out for a
fairly representative family of infinite \me-preserving \dsy s. These
are lattices of coupled baker's maps which generalize the \rw\ of
Section \ref{subs-ex} in two ways. First and foremost, they represent
\emph{all} the \rw s in $\Z^d$. Secondly, they are invertible \dsy s,
which can be reduced, for example, to the noninvertible \dsy\ of
Section \ref{subs-ex} by a mere restriction of the $\sigma$-algebra.

\skippar

To begin with, let a \rw\ in $\Z^d$ be defined by the transition
probabilities $\{ p_\b \}_{\b \in \Z^d}$, with $p_\b \ge 0$ and
$\sum_\b p_\b = 1$. This means that, if the walker is in $\a \in
\Z^d$, he will have \pr\ $p_\b$ to move to $\a+\b$ in the next
step. We introduce some notation that will be useful later:
\begin{equation}
  \label{def-dir}
  \D := \rset{\b \in \Z^d} {p_\b > 0} =: \{ \b^{(j)} \}_{j \in \Z_N}
\end{equation}
is the set of the ``active'' directions for the \rw, endowed with some
enumeration $\Z_N \ni j \mapsto \b^{(j)} \in \D$. If $\D$ is infinite,
then $N := \infty$ and $\Z_N := \Z^+$; if $\D$ is finite, then $N$
denotes its cardinality and $\Z_N := \Z^+ \cap [1,N]$.

We view this \rw\ as a \dsy\ $(\ps, \sca, \mu, T)$, where:
\begin{itemize}
  \item $\ps := \Z^d \times [0,1)^2$. If we denote $S_\a := \{\a\}
  \times [0,1)^2$, for $\a \in \Z^d$, then $\ps = \bigcup_\a S_\a$ can
  be interpreted as the disjoint union of $\Z^d$ copies of the unit
  square.

  \item $\sca$ is the natural $\sigma$-algebra for $\ps$, i.e., the
  $\sigma$-algebra generated by all the Lebesgue-measurable subsets of
  $S_\a$, $\forall \a \in \Z^d$, with the natural identification $S_\a
  \simeq [0,1)^2$.

  \item $\mu$ is the infinite \me\ that coincides with the Lebesgue
  \me\ when restricted to each $S_\a$.

  \item In order to define $T$, set $q_0 = 0$ and, for $k \in \Z_N$,
  \begin{equation}
    q_k := \sum_{j=1}^k p_{\b^{(j)}}; \qquad R_k := [q_{k-1}, q_k) 
    \times [0,1).
  \end{equation}
  Clearly $\{ R_k \}_{k \in \Z_N}$ is a partition of $[0,1)^2$ into
  adjacent rectangles of height 1 and width, respectively,
  \begin{equation}
    q_k - q_{k-1} = p_{\b^{(k)}} = \mu(R_k).
  \end{equation}
  For $x = (\a, y) = (\a; y_1, y_2) \in \Z^d \times [0,1)^2$, let $k
  \in \Z_N$ be the unique positive integer such that $y \in R_k$
  (equivalently, $q_{k-1} \le y_1 < q_k$). One defines
  \begin{equation}
    \label{def-map}
    Tx = T \left( \a; y_1, y_2 \right) := \left( \a + \b^{(k)} \,;
    \: p_{\b^{(k)}}^{-1} (y_1 - q_{k-1}) \,,\: p_{\b^{(k)}}   y_2 + 
    q_{k-1} \right).
  \end{equation}
\end{itemize}

Therefore $T$ is a piecewise linear, \hyp, invertible map $\ps
\longrightarrow \ps$ which preserves $\mu$ (because its determinant,
in the variables $(y_1, y_2)$, is 1). Denoting $R_{\a,k} := \{\a\}
\times R_k$, it is easy to see that $T$ is a Markov map for the
partition $\{ R_{\a,k} \}_{\a,k}$.

Now define $\psi: \ps \longrightarrow \Z^d$ as $\psi (\a,y) = \a$. It
is evident that, having chosen $x = (0,y)$ at random in $S_0$ w.r.t.\
$\mu$, the stochastic process $\{ \psi(T^n x) \}_{n\in \N}$ is the
\rw\ introduced at the beginning of the section, with initial position
in the origin.

\skippar

Moving on, we need to specify $\scv$, the exhaustive family of sets
that determines the infinite-volume limit: For all $\g = ( \g_1,
\ldots, \g_d ) \in \Z^d$ and $r \in \Z^+$, the set
\begin{equation}
  \label{def-box}
  B_{\g,r} = \rset{ \a = (\a_1, \ldots, \a_d ) \in \Z^d } {\g_i -r \le
  \a_i \le \g_i + r, \ \forall i=1,\ldots, d}
\end{equation}
is called a \emph{square box in $\Z^d$}. Then we pose
\begin{equation}
  \label{def-scv}
  \scv := \rset{ V = B_{\g,r} \times [0,1)^2} {\g \in \Z^d ,\ r \in
  \Z^+}
\end{equation}

\begin{remark}
  Clearly, definition (\ref{def-box}) does not capture all the square
  boxes in $\Z^d$, but only those whose side length is odd. This
  choice, made on grounds of simplicity, does not really limit the
  generality of $\scv$, and indeed the forthcoming results can be
  proven even in the case when (\ref{def-scv}) is modified to include
  all square boxes.
\end{remark}

\begin{lemma}
  \label{lem-a1}
  The \dsy\ and the exhaustive family defined above verify 
  {\bf (A1)}.
\end{lemma}

When $\D$ is finite, the proof of Lemma \ref{lem-a1} is rather
straightforward, along the same lines as (\ref{pf-a1-sec2}) for the
first example of Section \ref{subs-ex}. In the general case, it has
inessential technical complications, so we postpone it to Section
\ref{subs-pf-a1} of the Appendix.

\skippar

In order to introduce our \ob s, we need some preliminary
notation. Let $\scb \subset \sca$ be the $\sigma$-algebra on $\ps$
generated by the partition $\{ S_\a \}$ and, as is customary, $T^m
\scb = \rset{T^{-m} A} {A \in \scb}$. Then, for $\ell, m \in \Z$ with
$\ell \le m$, define
\begin{equation}
  \scb_{\ell,m} := T^\ell \scb \vee T^{\ell-1} \scb \cdots \vee
  T^{m+1} \scb \vee T^m \scb.
\end{equation}
To fix the ideas, $\scb_{0,1}$ is the $\sigma$-algebra corresponding
to the partition $\{ R_{\a,k} \}$. More generally, consider $\ell < 0
< m$: Recalling that $N = \# \D$ is the number of rectangles in each
partition $\{ R_k \}$ of $S_\a$, one can see that the fundamental
partition of $\scb_{\ell,m}$ is made up of $N^{|\ell| + m}$ rectangles
whose widths are bounded by $\lambda^m$ and whose heights are 
bounded by $\lambda^{|\ell|}$, where $\lambda := \max \{ p_\b \}$. 
Finally, set
\begin{eqnarray}
  \scb_{\ell, +\infty} &:=& \bigvee_{m \in \N} \scb_{\ell,m};
  \nonumber \\ 
  \scb_{-\infty, m} &:=& \bigvee_{-\ell \in \N} \scb_{\ell,m}; 
  \nonumber \\ 
  \scb_{-\infty, +\infty} &:=& \bigvee_{-\ell,m \in \N} \scb_{\ell,m}.
\end{eqnarray}
If we exclude, now and for the remainder of the section, the trivial
case where $N = 1$ (i.e., $p_\b = 0$, $\forall \b \ne \b^{(1)}$), one
has that $\lambda < 1$. This and the previous observation on the
fundamental sets of $\scb_{\ell,m}$ imply that the sets of
$\scb_{0,+\infty}$ are measurable unions of segments of the type $\{
\a \} \times \{ y_1 \} \times [0,1)$. We call those the \emph{local
stable manifolds} (LSMs) of the \sy\ and $\scb_{0,+\infty}$ the
\emph{stable $\sigma$-algebra}, also denoted by $\sca_s$. Analogously,
$\sca_u := \scb_{-\infty,0}$ is called the \emph{unstable
$\sigma$-algebra} and its sets are measurable unions of \emph{local
unstable manifolds} (LUMs) $\{ \a \} \times [0,1) \times \{ y_2 \}$.
Clearly, then, $\scb_{-\infty, +\infty} = \sca$.

We define several classes of global \ob s:
\begin{eqnarray}
  \label{def-gom}
  \go_m &:=& \rset{ F \in L^\infty (\ps, \scb_{-m,m}, \mu) } { \exists 
  \, \avg(F) \mbox{ as per definition \textbf{(A3)}} } ; \\
  \label{def-go}
  \go &:=& \wbar{\bigcup_{m \in \N} \go_m},
\end{eqnarray}
where the closure is meant in the $L^\infty$ norm. One should notice
that
\begin{lemma}
  \label{lem-go}
  Given the definitions \emph{(\ref{def-gom})-(\ref{def-go})},
  \begin{displaymath}
    \go = \lset{ F \in \wbar{ \bigcup_{m \in \N} L^\infty (\ps,
    \scb_{-m,m}, \mu) } } { \exists \, \avg(F) } .
  \end{displaymath}
\end{lemma}

\proof Let us prove the left-to-right inclusion. Given $F \in \go$ and
$n \ge 1$, there exists an $F_n \in \bigcup_m \go_m$ such that $\| F_n
- F \|_{L^\infty} \le 1/n$. This implies that
\begin{equation}
  \avg(F_n) - \frac1n \le \liminf_{V \nearrow \ps} \mu_V(F) \le 
  \limsup_{V \nearrow \ps} \mu_V(F) \le \avg(F_n) + \frac1n.
\end{equation}
On the other hand, $\{ \avg(F_n) \}$ is a Cauchy sequence, because $\{
F_n \}$ is. Its convergence thus proves the existence of $\avg(F)$.

Conversely, if $F$ is an \ob\ in the closure of $\bigcup_m L^\infty
(\ps, \scb_{-m,m}, \mu)$ for which $\avg(F)$ exists, then, for any
$\eps>0$, there exist $m \in \N$ and $F' \in L^\infty (\ps,
\scb_{-m,m}, \mu)$ such that
\begin{equation}
  \label{go-20}
  \| F' - F \|_{L^\infty} \le \eps/2. 
\end{equation}
Denoting $F'' := \mathbb{E} (F | \scb_{-m,m})$, (\ref{go-20}) implies
that $\| F' - F'' \|_{L^\infty} \le \eps/2$, hence $\| F - F''
\|_{L^\infty} \le \eps$. Furthermore, since $\scv \subset \scb \subset
\scb_{-m,m}$, then $\avg(F'')$ exists and equals $\avg(F)$. This shows
that $F'' \in \go_m$. Therefore $F \in \go$.
\qed

\begin{remark}
  In view of Lemma \ref{lem-go}, one might wonder why we do not
  consider the more natural class
  \begin{equation}
    \go_\infty := \lset{ F \in L^\infty (\ps, \sca, \mu) } 
    { \exists \, \avg(F) } 
  \end{equation}
  instead of $\go$. (The inclusion $\go \subset \go_\infty$ is clearly
  strict.) The reason, which will hardly surprise the ``hyperbolic''
  dynamicist, is that we need to approximate a global observable with
  locally constant \fn s \emph{uniformly} over $\ps$, cf.\
  (\ref{def-gom})-(\ref{def-go}). At any rate, $\go$ does not lack
  generality: for instance, any uniformly continuous $F$ verifying
  \textbf{(A3)} belongs in that class.
\end{remark}

As for the local \ob s, we also introduce countably many classes:
\begin{eqnarray}
  \label{def-lom}
  \lo_m &:=& L^1 (\ps, \scb_{-m,m}, \mu) ; \\
  \label{def-lo}
  \lo &:=& L^1 (\ps, \sca, \mu).
\end{eqnarray}

Prior to stating the main theorem of this section, we give a lemma
that will help appreciate its statement. If $\{ \a^{(j)} \}_{j \in
\mathcal{J}} \subset \Z^d$, the expression $\mathrm{span}_\Z \{
\a^{(j)} \}_{j \in \mathcal{J}}$ denotes the subgroup of all the
finite linear combinations of the $\a^{(j)}$ with coefficients in
$\Z$.

\begin{lemma}
  \label{lem-span}
  Let $\{ \b^{(j)} \}_{j \in \Z_N} \subset \Z^d$ and $j' \in \Z_N$.
  Then
  \begin{displaymath}
    \mathrm{span}_\Z \{ \b^{(j)} - \b^{(j')} \}_{j \in \Z_N}
  \end{displaymath}
  does not depend on $j'$.
\end{lemma}

\proof Section \ref{subs-pf-span} of the Appendix.

\begin{theorem}
  \label{thm-rw}
  Let $(\ps, \sca, \mu, T)$ be the \dsy\ described above. Set $\nu :=
  \max \{ 2, [d/2] + 1 \}$, where $[ \cdot ]$ the integer part of a
  positive number, and suppose that
  \begin{itemize}
    \item[(i)] the \pr\ distribution $p$ has a finite $\nu$-th momentum:
    \begin{displaymath}
      \sum_{\b \in \Z^d} |\b|^\nu p_\b < \infty ;
    \end{displaymath}

    \item[(ii)] for a given $j' \in \Z_N$, $\mathrm{span}_\Z \{ \b^{(j)}
    - \b^{(j')} \}_{j \in \Z_N} = \Z^d$ .
  \end{itemize}
  Then the \sy\ is mixing in the following senses:
  \begin{itemize}
    \item[(a)] \emph{\textbf{(M5)}} relative to $\go_m$ and $\lo_m$, for
    all $m \in \N$;

    \item[(b)] \emph{\textbf{(M4)-(M3)}} relative to $\go$ and $\lo$;

    \item[(c)] \emph{\textbf{(M2)}} relative to $\go$;

    \item[(d)] \emph{\textbf{(M1)}} relative to $\go_m$, for all $m
    \in \N$, with the extra requirement that $F$ be $\Z^d$-periodic,
    i.e., $F(\a,y) = F(0,y)$, $\forall \a \in \Z^d$, $\forall y \in
    [0,1)^2$.
  \end{itemize}
\end{theorem}

\begin{remark}
  Condition \emph{(ii)} is essential as it has to do with the
  irreducibility of the \rw\ \cite{sp}. In fact, assuming for
  simplicity that $0 \in \D$, if $\mathrm{span}_\Z (\D) \ne \Z^d$, then
  the \rw\ is reducible and the \sy\ cannot be mixing in any sense, as
  is ascertained via the global observable
  \begin{equation}
    F(\a, y) := \left\{
    \begin{array}{lll}
      1, && \mbox{for } \a \in \mathrm{span}_\Z (\D) ; \\
      0, && \mbox{otherwise}.
    \end{array}
  \right.
  \end{equation}
\end{remark}

One last observation that may be of interest is that statement
\emph{(d)} is far from optimal. \textbf{(M1)} holds for a much larger
class of global \ob s, depending especially on the distribution $\{
p_\b \}$. In the formulation of Theorem \ref{thm-rw}, however, I was
mainly interested in a nontrivial case in which \textbf{(M1)} could be
verified easily.

\section{Proof of Theorem \ref{thm-rw}}
\label{sec-pf-main}

Since the proof is rather lengthy, we will divide it into pieces, or
\emph{stages}, as follows:

\smallskip\skippar
\emph{Stage 1:} We prove \emph{(a)} using three extra assumptions.
\skippar
\emph{Stage 2:} We remove one of the extra assumptions.
\skippar
\emph{Stage 3:} We remove the remaning two extra assumptions.
\skippar
\emph{Stage 4:} We prove \emph{(b)-(d)}.

\subsection{Stage 1: Extra assumptions}

Let us initially assume that:

\begin{itemize}
\item[(E1)] $F \in \sca_s = \scb_{0, +\infty}$;

\item[(E2)] $\lo_m$ only comprises indicator \fn s of the type $g =
  1_Q$, where $Q$ is a fundamental set of $\scb_{-m,0}$ and $Q
  \subseteq S_0$;

\item[(E3)] the \rw\ has zero drift, i.e., $\ds \sum_{\b \in
  \Z^d} \b \, p_\b = 0$.
\end{itemize}

From (E2), $Q$ is a rectangle of the type $\{ 0 \} \times
[0,1) \times I$. We denote the length of $I$ by
\begin{equation}
  \label{t-rw-5}
  h := |I| = \mu(Q).
\end{equation}
In this setting, \textbf{(M5)} amounts to showing that
\begin{equation}
  \label{t-rw-goal}
  \lim_{n \to \infty} \, \frac1 {\mu(|g|)} \int_\ps (F \circ T^n) g \, d\mu
  = \lim_{n \to \infty} \, \frac1h \int_{T^n Q} F \, d\mu = \avg(F),
\end{equation}
\emph{uniformly in $Q$}, that is, with a speed of convergence that
does not depend on the choice of $I$. (In the above we have used the
invariance of $\mu$ and the fact that $\mu(|g|) = h = \mu(g)$.)
Achieving this will be Stage 1 of the proof.

\skippar

$Q$ can be thought of as partitioned into $\{ Q \cap R_{0,j} \}_{j \in
\Z_N}$, which are rectangles of width $p_{\b^{(j)}}$ and height
$h$. By construction, $T$ acts on each such rectangle by stretching it
horizontally and shrinking it vertically by a factor
$p_{\b^{(j)}}^{-1}$, and then mapping the resulting rectangle, of
width 1, rigidly into $S_{\b^{(j)}}$.

Iterating this procedure $n$ times, we obtain
\begin{equation}
  \label{t-rw-10}
  T^n Q = \bigcup_{j_1, \ldots, j_n \in \Z_N} Q_{j_1, \ldots, j_n}
  := \bigcup_{j_1, \ldots, j_n \in \Z_N} \{ \a_{j_1, \ldots, j_n} \}
  \times [0,1) \times I_{j_1, \ldots, j_n} ,
\end{equation}
which is a disjoint union of $N^n$ thin rectangles of width 1.  Each
$Q_{j_1, \ldots, j_n}$ is the set of all the points $T^n x$ for which
the \tr y of $x \in Q$ followed the \emph{itinerary} $S_{ \b^{(j_1)}
}$, $S_{ \b^{(j_1)} + \b^{(j_2)} }$, \ldots, $S_{ \a_{j_1, \ldots,
    j_n} }$, where
\begin{equation}
  \label{t-rw-30}
  \a_{j_1, \ldots, j_n} := \b^{(j_i)} + \b^{(j_2)} + \cdots +
  \b^{(j_n)}.
\end{equation}
Therefore, recalling (\ref{t-rw-5}), the height of $Q_{j_1, \ldots,
j_n}$ is 
\begin{equation}
  \label{t-rw-40}
  h_{j_1, \ldots, j_n} := |I_{j_1, \ldots, j_n}| = \mu( Q_{j_1,
  \ldots, j_n} ) = h \, p_{\b^{(j_1)}} p_{\b^{(j_2)}} \cdots \,
  p_{\b^{(j_n)}}.
\end{equation}
Since $F \in \sca_s$, we can write, with a slight abuse of notation,
$F(\a; y_1, y_2) = F_\a(y_1)$. Then
\begin{eqnarray}
  \int_{T^n Q} F \, d\mu &=& \sum_{j_1, \ldots, j_n \in \Z_N} \:
  \int_{Q_{j_1, \ldots, j_n}} F \, d\mu = \nonumber \\ 
  &=& \sum_{j_1, \ldots, j_n \in \Z_N} \int_{I_{j_1, \ldots, j_n}} \:
  \int_0^1 F_{\a_{j_1, \ldots, j_n}} (y_1) \, dy_1 \ dy_2 = \nonumber
  \\
  &=& \sum_{j_1, \ldots, j_n \in \Z_N} h \, p_{\b^{(j_1)}} p_{\b^{(j_2)}}
  \cdots \, p_{\b^{(j_n)}} \, f_{\a_{j_1, \ldots, j_n}},
  \label{t-rw-50}
\end{eqnarray}
having used (\ref{t-rw-10}), (\ref{t-rw-40}) and the following
definition:
\begin{equation}
  \label{t-rw-60}
  f_\a := \int_0^1 F_\a (y_1) \, dy_1.
\end{equation}
In view of (\ref{t-rw-30}) and (\ref{def-dir}), another way to write
(\ref{t-rw-50}) is
\begin{equation}
  \label{t-rw-70}
  \int_{T^n Q} F \, d\mu = h \sum_{\b^{(1)}, \ldots, \b^{(n)} \in
  \Z^d} p_{\b^{(1)}} p_{\b^{(2)}} \cdots \, p_{\b^{(n)}} \,
  f_{\b^{(1)} + \ldots + \b^{(n)}}.
\end{equation}

\subsection{Fourier analysis}

The technical backbone of the proof is Fourier analysis on $\Z^d$, for
which we proceed to establish the necessary notation \cite{k, r}. Let
$a := \{ a_\a \}_{\a \in \Z^d} \in \ell^s(\Z^d; \C)$, with $s \in
[1,\infty]$. Its Fourier transform is denoted
\begin{equation}
  \label{t-rw-ft}
  \wtilde{a}(\t) = \wtilde{a} (\t_1, \ldots, \t_d) := \sum_{\a \in
  \Z^d} a_\a \, \ei{ \a \cdot \t } = \sum_{\a \in \Z^d} a_\a \, \ei{ 
  ( \a_1 \t_1 + \ldots + \a_d \t_d ) },
\end{equation}
where $\imath$ is the imaginary unit and $\t = ( \t_1, \ldots, \t_d )
\in \T^d := (\R / 2\pi\Z)^d$. If $s>2$, (\ref{t-rw-ft}) must be
intended in the weak sense, cf.~(\ref{t-rw-pars}). The corresponding
inverse transform is given by
\begin{equation}
  \label{t-rw-ift}
  a_\a = \int_{\T^d} \wtilde{a}(\t) \, \emi{ \a \cdot \t } \,
  \dpt,
\end{equation}
where $\dpt := (2\pi)^{-d} d\t$. The Parseval formula, in this
setting, reads as follows: If $b = \{ b_\a \} \in \ell^{s'} (\Z^d;
\C)$, with $1/s + 1/s' = 1$, then
\begin{equation}
  \label{t-rw-pars}
  \langle a,b \rangle := \sum_{\a \in \Z^d} \wbar{a_\a} \, b_\a =
  \int_{\T^d} \wbar{ \wtilde{a}(\t) } \: \wtilde{b}(\t) \, \dpt,
\end{equation}
the bar denoting complex conjugation. Another standard result that we
need is the duality between convolution and product: If
\begin{equation}
  \label{t-rw-defconv}
  (a * b)_\a := \sum_{\b \in \Z^d} a_\b \, b_{\a-\b} = \sum_{\b \in
  \Z^d} a_{\a-\b} \, b_\b
\end{equation}
is well-defined in the proper or weak sense, then
\begin{equation}
  \label{t-rw-ftconv}
  \wtilde{(a * b)} (\t) = \wtilde{a} (\t) \, \wtilde{b} (\t).
\end{equation}

Applying these concepts to our case, we see that $f := \{ f_\a \} \in
\ell^\infty$ by construction (because $F \in \go \subset
L^\infty(\ps)$) so $\wtilde{f}$ is a distribution on $\T^d$. On the
other hand, $p := \{ p_\a \} \in \ell^1$, which makes $\wtilde{p}$
continuous (it is actually much more than that, cf.~Section
\ref{subs-prop}). In particular, since $p$ is a \pr\ distribution on
$\Z^d$, $\wtilde{p}(0) = 1$. Defining
\begin{equation}
  \label{t-rw-80}
  p^{(n)} := \underbrace{p * \cdots * p}_{n\ \mathrm{times}} 
\end{equation}
we can rearrange (\ref{t-rw-70}) into
\begin{equation}
  \label{t-rw-90}
  \frac1{h} \int_{T^{-n} Q} F \, d\mu = \left\langle f , p^{(n)}
  \right\rangle = \int_{\T^d} \wbar{ \wtilde{f} (\t) } \, \wtilde{
  p^{(n)} } (\t) \, \dpt = \int_{\T^d} \wbar{ \wtilde{f} (\t) } \,
  \wtilde{p}^n (\t) \, \dpt,
\end{equation}
where we have used (\ref{t-rw-ftconv}), (\ref{t-rw-80}) and the fact
that $f_\a \in \R$. Hence Stage 1 of the proof, cf.~(\ref{t-rw-goal}),
reduces to showing that
\begin{equation}
  \label{t-rw-goal2}
  \lim_{n \to \infty} \left\langle f , p^{(n)} \right\rangle =
  \avg(F).
\end{equation}
Now recall definition (\ref{def-box}). For $r \in \N$, let 
\begin{equation}
  \label{t-rw-100}
  q^{(r)}_\a := \left\{
  \begin{array}{lll}
    (2r+1)^{-d}, && \mbox{if } \a \in B_{0,r} ; \\
    0, && \mbox{otherwise} ,
  \end{array} \right.
\end{equation}
define a \fn\ $q^{(r)} : \Z^d \longrightarrow \R$. Its Fourier
transform is easily computed to be
\begin{equation}
  \label{t-rw-110}
  \wtilde{ q^{(r)} }( \t_1, \ldots, \t_d ) = \prod_{i=1}^d \frac{ 
  \sin \left( (r + 1/2) \t_i \right) } { (r + 1/2) \sin \t_i } .
\end{equation}

In view of (\ref{def-scv}), let us denote $V_{\a,r} := B_{\a,r} \times
[0,1)^2$. Since $F$ verifies \textbf{(A3)} and the infinite-volume
limit is $\mu$-uniform (cf.\ Definition \ref{def-iv}), we have that
\begin{eqnarray}
  \label{t-rw-120}
  && \lim_{r \to \infty} \frac1 {\mu(V_{\a,r})} \int_{V_{\a,r}} \! F 
  \, d\mu = \lim_{r \to \infty} \frac1{(2r+1)^d} 
  \sum_{\b \in B_{\a,r}} f_\b = \nonumber \\
  && \lim_{r \to \infty} \left( f * q^{(r)} \right)_\a = \avg(F),
\end{eqnarray}
\emph{uniformly in $\a$}. (We have used notation (\ref{t-rw-60}).)
This, in turn, yields
\begin{equation}
  \label{t-rw-130}
  \lim_{r \to \infty} \left\langle f , q^{(r)} * p^{(n)}
  \right\rangle = \lim_{r \to \infty} \left\langle f * q^{(r)} ,
  p^{(n)} \right\rangle = \avg(F)
\end{equation}
\emph{uniformly in $n$}, because $p^{(n)}$ is a \pr\ distribution on
$\Z^d$. (In the first equality we have used the fact that $q^{(r)}_\a
= q^{(r)}_{-\a}$, by construction.) This fact implies that, for any
sequence $\{ r_n \} \subset \N$ with $r_n \to \infty$,
\begin{equation}
  \label{t-rw-132}
  \lim_{n \to \infty} \left\langle f , q^{(r_n)} * p^{(n)}
  \right\rangle =\avg(F).
\end{equation}
Therefore, comparing (\ref{t-rw-132}) with (\ref{t-rw-goal2}), we see
that Stage 1 is achieved once we have shown that there exists a
diverging sequence $\{ r_n \}$ of natural numbers such that
\begin{equation}
  \label{t-rw-goal3}
  \lim_{n \to \infty} \left\langle f , p^{(n)} - q^{(r_n)} * p^{(n)}
  \right\rangle = 0.
\end{equation}

For an \emph{a fortiori} choice of $\{ r_n \}$, let us set
\begin{equation}
  \label{t-rw-135}
  g^{(n)}_\a := p^{(n)}_\a - \left( q^{(r_n)} * p^{(n)} \right)_\a,
\end{equation}
which gives
\begin{equation}
  \label{t-rw-140}
  \wtilde{g^{(n)}} (\t) = \left( 1 - \wtilde{ q^{(r_n)} } (\t) \right) 
  \wtilde{p}^n (\t).
\end{equation}
A convenient estimate in view of (\ref{t-rw-goal3}) is 
\begin{equation}
  \label{t-rw-150}
  \left| \left\langle f , g^{(n)} \right\rangle \right| \le \left\| f
  \right\|_{\ell^\infty} \left\| g^{(n)} \right\|_{\ell^1} = \left\| f
  \right\|_{\ell^\infty} \left\| \wtilde{g^{(n)}} \right\|_{\ac}
  \le C \left\| \wtilde{ g^{(n)} } \right\|_{H^\nu},
\end{equation}
where the norms $\| \,\cdot\, \|_{\ac}$ and $\| \,\cdot\, \|_{H^\nu}$
are introduced in Section \ref{subs-anorm} of the Appendix (cf.\ in
particular Lemma \ref{lem-e-anorm} and notice that $\| \,\cdot\,
\|_{H^{\bar{\nu}}} \le \| \,\cdot\, \|_{H^\nu}$).

Therefore (\ref{t-rw-goal3}) will be proved once we establish that
\begin{equation}
  \label{t-rw-goal4}
  \left\| \wtilde{g^{(n)}} \right\|_{H^\nu} \le \left\|
  \wtilde{g^{(n)}} \right\|_{L^1} + \sum_{i=1}^d \left\|
  \partial_i^\nu \, \wtilde{g^{(n)}} \right\|_{L^2} \to 0, \quad
  \mbox{as } n \to \infty,
\end{equation}
for a suitable choice of $\{ r_n \}$ in definition (\ref{t-rw-135}).
Here $\partial_i := \partial / \partial\t_i$, acting on \fn
s $\T^d \longrightarrow \C$.

\subsection{Properties of $\wtilde{ q^{(r)} }$ and $\wtilde{p}$}
\label{subs-prop}

In view of the above goal we need to study some properties of the \fn
s $\wtilde{ q^{(r)} }$ and $\wtilde{p}$.

\begin{remark}
  \label{rk-no-ea3}
  None of the proofs in this section will use (E3).
\end{remark}

As a start, let us notice that $\wtilde{ q^{(r)} }$ is $C^\infty$ by
construction and $\wtilde{p}$ is $C^\nu$, with $\nu \ge 2$, by hypothesis
\emph{(i)} of Theorem \ref{thm-rw}.

\begin{lemma}
  \label{lem-pq2}
  Fix $r \in \Z^+$. On $\T^d$, $\wtilde{ q^{(r)} } (0) = \wtilde{p}(0)
  = 1$ and, for $\t \ne 0$,
  \begin{displaymath}
    |\wtilde{ q^{(r)} } (\t)| <1, \qquad |\wtilde{p}(\t)| < 1.
  \end{displaymath}
\end{lemma}

\proof We first prove the assertions on $\wtilde{p}$. Since $p$ is a
\pr\ distribution, $\wtilde{p}(0) = 1$ and $|\wtilde{p}(\t)| \le 1$,
$\forall \t$. Suppose by contradiction that $\exists \t' \in \T^d$,
$\t' \ne 0$, such that $|\wtilde{p}(\t')| = 1$, that is,
\begin{equation}
  \label{pq2-10}
  \wtilde{p}(\t') = \sum_{j \in \Z_N} p_{\b^{(j)}} \, \ei{\t' \cdot
  \b^{(j)}} = \ei{a},
\end{equation}
for some $a \in \R$. Since $p_{\b^{(j)}} > 0$ and $\sum_j p_{\b^{(j)}}
= 1$, necessarily $\ei{\t' \cdot \b^{(j)}} = \ei{a}$, $\forall j \in
\Z_N$, whence, $\forall j, j'$,
\begin{equation}
  \label{pq2-20}
  \ei{ \t' \cdot (\b^{(j)} - \b^{(j')}) } = 1.
\end{equation}

Let us define the \emph{character} (i.e., the homomorphism $\Z^d
\longrightarrow S^1 \subset \C$) $\eta_{\t'} (\a) := \ei{\t' \cdot
\a}$. It is easy to see that $\eta_{\t'}$ is not the trivial
character (which instead corresponds to $\t' = 0$; this is a
particular case of the so-called Pontryagin Duality
\cite[Thm.~2.1.2]{r}). On the other hand, (\ref{pq2-20}) reads
$\eta_{\t'} (\b^{(j)} - \b^{(j')}) = 1$ and hypothesis \emph{(ii)} of
Theorem \ref{thm-rw} implies that $\eta_{\t'} \equiv 1$, thereby
creating a contradiction.

As for the assertions on $\wtilde{ q^{(r)} }$, there is nothing more
to prove, because $\wtilde{ q^{(r)} }$ satisfies the same properties
as $\wtilde{p}$, insofar as the above argument is concerned.
\qed

\noindent
\textbf{Notational convention.} \textit{From now on, $\con$ will
  denote a generic universal constant. This means that its actual
  value will vary from formula to formula but will never depend on
  $n$, $r$, or $\t$.}

\begin{lemma}
  \label{lem-p3}
  If we think of $\wtilde{p}$ as a periodic \fn\ on $\R^d$ (as opposed
  to a \fn\ on $\T^d$), there exists a neighborhood $\mathcal{U}$ of
  $\t = 0$ and a positive constant $\con$ such that, for $\t \in
  \mathcal{U}$,
  \begin{displaymath}
    | \wtilde{p}(\t) | \le 1 - \con |\t|^2 = 1 - \con \left( \t_1^2 + 
    \cdots + \t_d^2 \right).
  \end{displaymath}
\end{lemma}

\proof As $\t \to 0$,
\begin{equation}
  \label{p3-10}
  \wtilde{p}(\t) = 1 + \imath v \cdot \t + \mathcal{O}(|\t|^2),
\end{equation}
where $v \in \R^d$ is the \emph{drift} of the \rw, defined as
\begin{equation}
  \label{p3-20}
  v := \sum_{\b \in \Z^d} \b p_\b. 
\end{equation}
The Lagrange remainder in (\ref{p3-10}) holds because $\wtilde{p}$ is
at least $C^2$. Hence $| \wtilde{p}(\t) |^2 = 1 +
\mathcal{O}(|\t|^2)$, which implies the assertion.
\qed

\begin{lemma}
  \label{lem-q3}
  Regarding $\wtilde{ q^{(r)} }$ as a periodic \fn\ on $\R^d$, one has
  that the following expansion, 
  \begin{displaymath}
    \wtilde{ q^{(r)} }( \t_1, \ldots, \t_d ) = \prod_{i=1}^d \left( 
    1 + \sum_{j=1}^\infty \xi_j(r) \, \t_i^{2j} \right),
  \end{displaymath}
  holds uniformly on the compact subsets of $\R^d$. Furthermore, as $r
  \to \infty$, 
  \begin{displaymath}
    |\xi_j(r)| \le C \frac{r^{2j}} {(2j)!}.
  \end{displaymath}
\end{lemma}

\proof By the factorizability of $q^{(r)}$, it is sufficient to treat
the case $d = 1$.

Since $q^{(r)}$ is compactly supported in $\Z^d$, $\wtilde{ q^{(r)} }
(\t)$ is an entire \fn\ of $\t$, which we have already calculated in
(\ref{t-rw-100}). Its Taylor expansion at the origin is even and its
(even) terms are
\begin{equation}
  \label{q3-10}
  \xi_j(r) = \frac1 {(2j)!} \, \partial^{2j} \wtilde{q^{(r)}} (0) =
  \frac1 {(2j)!} \, \frac{(-1)^j} {2r+1} \sum_{\a = -r}^r \a^{2j}.
\end{equation}
This gives $\xi_0(r) \equiv 1$ and the desired estimates.
\qed

We estimate the norms in the r.h.s.\ of (\ref{t-rw-goal4}) by
splitting the corresponding integrals into two parts: one over
$\mathcal{B}_n$, which is the ball of center $0$ and radius
$n^{-(1-\eps)/2}$ in $\T^d$, and the other over $\T^d \setminus
\mathcal{B}_n$. Here $\eps > 0$ is a small constant to be fixed later
and $n$ is a large integer.

\begin{lemma}
  \label{lem-p4}
  There exists a $\kappa > 0$ such that, for $n$ sufficiently large,
  \begin{displaymath}
    \max_{\T^d \setminus \mathcal{B}_n} |\wtilde{p}|^n \le 
    e^{-\kappa n^\eps}.
  \end{displaymath}
\end{lemma}

\proof By elementary Taylor approximations, Lemma \ref{lem-p3} implies
that there exists a constant $\kappa>0$ such that, for $\t \in
\mathcal{U}$, 
\begin{equation}
  \label{p4-10}
  | \wtilde{p}(\t) | \le e^{-\kappa |\t|^2}.
\end{equation}
For $n$ large enough, by Lemma \ref{lem-pq2}, the continuity of
$\wtilde{p}$ and the compactness of $\T^d$,
\begin{equation}
  \label{p4-20}
  \max_{\t \in \T^d \setminus \mathcal{B}_n} | \wtilde{p} (\t) | = 
  \max_{\t \in \partial \mathcal{B}_n} | \wtilde{p} (\t) | \le
  e^{-\kappa n^{-1+\eps}}
\end{equation}
(in the last inequality we have used (\ref{p4-10}), which applies
because, for $n$ large, $\mathcal{B}_n \subset \mathcal{U}$).
\qed

\subsection{End of Stage 1}
\label{subs-endst1}

Let us begin to attack (\ref{t-rw-goal4}) by estimating $\|
\wtilde{g^{(n)}} \|_{L^1}$. First of all,
\begin{equation}
  \label{t-rw-160}
  \left\| \wtilde{g^{(n)}} \right\|_{ L^1 (\T^d \setminus
  \mathcal{B}_n) } \le (2\pi)^d 2 e^{-\kappa n^\eps},
\end{equation}
by Lemma \ref{lem-pq2} applied to $\wtilde{q^{(r)}}$ and Lemma
\ref{lem-p4}---see (\ref{t-rw-140}). Again, applying Lemma
\ref{lem-pq2} to both $\wtilde{q^{(r)}}$ and $\wtilde{p}$,
\begin{equation}
  \label{t-rw-170}
  \left\| \wtilde{g^{(n)}} \right\|_{ L^1 (\mathcal{B}_n) } \le
  \frac{\con} { n^{d(1-\eps)/2} } 
\end{equation}
where $C$ does not depend on $r$, that is, $r_n$. 

As for the remaining terms in the r.h.s.\ of (\ref{t-rw-goal4}),
clearly
\begin{equation}
  \label{t-rw-180}
  \partial_i^\nu \wtilde{g^{(n)}} = \sum_{k+l = \nu} {\nu \choose k}
  \: \partial_i^k \! \left( 1 - \wtilde{q^{(r)}} \right) \partial_i^l
  \, \wtilde{p}^n.
\end{equation}
Let us estimate (\ref{t-rw-180}) on $\T^d \setminus \mathcal{B}_n$.
Fixing $l \ge 1$, one verifies by repeated differentiation that
\begin{equation}
  \label{t-rw-190}
  \partial_i^l \, \wtilde{p}^n = \sum_{w=1}^l \, \sum_{ j_1 + \ldots +
  j_w = l } C_{j_1, \ldots, j_w}^l \, \frac{n!} {(n-w)!} \,
  \wtilde{p}^{n-w} (\partial_i^{j_1} \wtilde{p}) (\partial_i^{j_2}
  \wtilde{p}) \cdots (\partial_i^{j_w} \wtilde{p}),
\end{equation}
where the combination of the two sums above represents the sum over
all the partitions $\{ j_1, j_2, \ldots, j_w \}$ of $l$ (i.e., $j_u
\ge 1$ and $j_1 + j_2 + \ldots + j_w = l$), with any cardinality $w$,
and $C_{j_1, \ldots, j_w}^l \in \N$ is a combinatorial coefficient
independent of $n$. Since $\wtilde{p} \in C^\nu$, all the derivatives
that appear on the r.h.s.\ of (\ref{t-rw-190}) are continuous \fn s of
$\t$. Therefore, by Lemma \ref{lem-p4},
\begin{equation}
  \label{t-rw-193}
  \max_{\T^d \setminus \mathcal{B}_n} \left| \partial_i^l \,
  \wtilde{p}^n \right| \le \con n^l e^{-\kappa (n - l)^\eps}.
\end{equation}
As for the other factors in the r.h.s.\ of (\ref{t-rw-180}), for $k
\ge 1$ we use definition (\ref{t-rw-100}) to estimate
\begin{equation}
  \label{t-rw-197} 
  \max_{\T^d} \left| \partial_i^k \wtilde{q^{(r)}} \right| \le
  \sum_{\a \in \Z^d} \left| \a_i^k q_\a^{(r)} \right| = \frac1 {2r+1}
  \sum_{\a_i = -r}^r |\a_i|^k \le \con r^k.
\end{equation}
(For $k=0$, we already know that $1 - \wtilde{q^{(r)}}$ is bounded.)
Using (\ref{t-rw-193})-(\ref{t-rw-197}) into (\ref{t-rw-180}), we get 
\begin{equation}
  \label{t-rw-200}
  \left\| \partial_i^\nu \wtilde{g^{(n)}} \right\|_{ L^2 (\T^d
  \setminus \mathcal{B}_n) }^2 \le \con r^{2\nu} n^{2\nu} e^{2\kappa
  (n - \nu)^\eps},
\end{equation}
which tends to zero, as $n \to \infty$, provided that $r = r_n$ grows
no faster than a power of $n$. This will be verified \emph{a
fortiori}, see (\ref{t-rw-rn}).

The estimation of the last term, namely $\| \partial_i^\nu
\wtilde{g^{(n)}} \|_{L^2 (\mathcal{B}_n)}$, is the most delicate,
therefore we organize most of the computations involved in the
following

\begin{lemma}
  \label{lem-last}
  For $\nu \in \Z^+$ (not necessarily as in the statement of Theorem 
  \ref{thm-rw}), assume $\wtilde{p} \in C^\nu$. Then take any
  sequence of positive numbers $\Lambda_n$, with $\Lambda_n \to 0$. For
  $n$ large enough and uniformly for
  \begin{displaymath}
    1 \le r \le \Lambda_n \, n^{(1-\eps)/2},
  \end{displaymath}
  one has
  \begin{displaymath}
    \max_{\mathcal{B}_n} \left| \partial_i^\nu \wtilde{g^{(n)}}
    \right| \le \con r^{2\nu} n^{(-2 + 2\eps + \nu (1+\eps))/2}.
  \end{displaymath}
\end{lemma}

Once Lemma \ref{lem-last} is established, which will happen
momentarily, we can finally choose both the sequence $\{ r_n \}$ and
the parameter $\eps$: the former will be any diverging sequence such
that
\begin{equation}
  \label{t-rw-rn}
  r_n \le \con n^\eps,
\end{equation}
whereas 
\begin{equation}
  \label{t-rw-eps}
  \eps < \min \left\{ \frac13 ,\, \frac1 {2(5\nu + 2 + d/2)} \right\}.
\end{equation}
The condition $\eps < 1/3$ is necessary in order to apply Lemma
\ref{lem-last} to $r = r_n$. In fact, via (\ref{t-rw-rn}), $r_n
n^{-(1-\eps)/2} \le \con n^{-(1-3\eps)/2} =: \Lambda_n$. The other
inequality is needed when we use the assertion of Lemma \ref{lem-last}
to estimate
\begin{eqnarray}
  \left\| \partial_i^\nu \wtilde{g^{(n)}} \right\|_{ L^2
  (\mathcal{B}_n) }^2 &\le& \con r_n^{4\nu} \, n^{-2 + 2\eps + \nu
  (1+\eps)} \, n^{-d (1-\eps) /2} \le \nonumber \\
  \label{t-rw-210}
  &\le& \con n^{-2 + \nu - d/2 + (5\nu + 2 + d/2)\eps}.
\end{eqnarray}
For $\eps$ so small that $(5\nu + 2 + d/2) \eps < 1/2$ the above
vanishes, as $n \to \infty$, because $\nu - d/2 \le 3/2$ (for $d = 1$,
it equals $3/2$; for $d \ge 2$ it equals $1$ or $1/2$, depending on
$d$ being even or odd).

In conclusion, estimates (\ref{t-rw-160}), (\ref{t-rw-170}),
(\ref{t-rw-200}) and (\ref{t-rw-210}) prove (\ref{t-rw-goal4}), thus
(\ref{t-rw-goal3}), thus (\ref{t-rw-goal2}) and, lastly,
(\ref{t-rw-goal}), uniformly in $Q$ as requested. This ends Stage 1 of
the proof.

\skippar

\proofof{Lemma \ref{lem-last}} Throughout the proof it is understood
that $\t \in \mathcal{B}_n$, i.e., $|\t| \le n^{-(1-\eps)/2}$. The
condition on $r$ ensures that
\begin{equation}
  \label{t-rw-220}
  \left| \xi_j(r) \t^{2j} \right| \le \con r^{2j} n^{-(1-\eps)j}
  \le \con \Lambda_n^{2j},
\end{equation}
which implies that, for $n$ so large that $\Lambda_n < 1$, the
expansion of Lemma \ref{lem-q3} and all those derived by it by
differentiation w.r.t.\ $\t_i$ are meaningful. More importantly, as $n
\to \infty$ and uniformly in $r$ as described in the statement of the
lemma, each such expansion can be approximated by a convenient upper
bound on its first term. We indicate this with the symbol $\sim$: for
example,
\begin{equation}
  \label{t-rw-230}
  1 - \wtilde{q^{(r)}} (\t) \,\sim\, \left| \xi_1(r) \t \right|^2
  \le \con r^2  n^{-(1-\eps)}
\end{equation}
and, for $k \ge 1$,
\begin{equation} 
  \label{t-rw-240}
  \left| \partial_i^k \wtilde{q^{(r)}} (\t) \right| \,\sim\, \left\{
  \begin{array}{lcll}
    (k+1)! \: | \xi_{(k+1)/2}(r) \t_i| &\le& \con r^{k+1}
    n^{-(1-\eps)/2}, & \mbox{ if $k$ is odd,} \\
    k! \: | \xi_{k/2}(r) | &\le& \con  r^k, & \mbox{ if $k$ is even.}
    \\ 
  \end{array}
  \right. 
\end{equation}
Furthermore, (E3) is equivalent to $\nabla \wtilde{p} (0) =
0$, which implies
\begin{equation}
  \label{t-rw-250}
  | \partial_i \wtilde{p} (\t) | \le \con |\t| \le \con
  n^{-(1-\eps)/2}.
\end{equation}

We will proceed by induction on $\nu \ge 1$. When $\nu = 1$ our \fn\
reads
\begin{equation}
  \label{t-rw-260}
  \partial_i \wtilde{g^{(n)}} = -\partial_i \wtilde{q^{(r)}} \,
  \wtilde{p}^n + \left( 1 - \wtilde{q^{(r)}} \right) n
  \wtilde{p}^{n-1} \, \partial_i \wtilde{p}.
\end{equation}
Hence, from (\ref{t-rw-230})-(\ref{t-rw-250}), and Lemma
\ref{lem-pq2} applied to $\wtilde{p}$,
\begin{equation}
  \label{t-rw-270}
  \max_{\mathcal{B}_n} \left| \partial_i \wtilde{g^{(n)}} \right| \le
  \con r^2 n^{(-1 +3\eps)/2},
\end{equation}
proving the assertion for $\nu=1$.

Now we assume the assertion with $\nu$ and set out to prove the one
with $\nu+1$. In practice, this means that increasing the order of the
derivative by one must worsen the inequality of Lemma \ref{lem-last}
at most by a factor $r^2 n^{(1+\eps)/2}$.

We apply $\partial_i$ to (\ref{t-rw-180}). On the r.h.s., $\partial_i$
can either hit $\partial_i^k ( 1 - \wtilde{q^{(r)}} )$ or
$\partial_i^l \wtilde{p}^n$. Let us analyze the two cases separately.

In the first case, assuming for the moment $k \ge 1$, we see via
(\ref{t-rw-240}) that
\begin{equation} 
  \label{t-rw-280}
  \left| \partial_i^{k+1} \wtilde{q^{(r)}} (\t) \right| \,\le\,
  \left\{
  \begin{array}{ll}
    Cr^{k+1} = \left( \con r^{k+1} \, n^{-(1-\eps)/2} \right)
    n^{(1-\eps)/2}, & \mbox{ if $k$ is odd,} \\
    Cr^{k+2} n^{-(1-\eps)/2} = \left( \con r^k \right) r^2 \,
    n^{-(1-\eps)/2}, & \mbox{ if $k$ is even.} \\ 
  \end{array}
  \right. 
\end{equation}
The terms within parentheses in the above represent the estimates
for $\partial_i^k \wtilde{q^{(r)}}$, respectively for $k$ odd and
even, coming from (\ref{t-rw-240}). Also, for $k=0$,
\begin{equation}
  \label{t-rw-300}
  \partial_i \wtilde{q^{(r)}} \le \con r^2 \, n^{-(1-\eps)/2} = 
  \left( \con r^2 \, n^{-(1-\eps)} \right) n^{(1-\eps)/2}.
\end{equation}
Again, the term in the parentheses is the estimate for $1 -
\wtilde{q^{(r)}}$ coming from (\ref{t-rw-230}).  In any event,
applying another derivative to the term $\partial_i^k ( 1 -
\wtilde{q^{(r)}} )$ will change our estimate at most by a factor $r^2
n^{(1 - \eps)/2}$, which is consistent with our inductive step.

In the second case, we use the expansion (\ref{t-rw-190}): 
$\partial_i$ can either hit $\wtilde{p}^{n-w}$ or one of the
$\partial_i^{j_u} \wtilde{p}$. In this first sub-case,
\begin{equation}
  \label{t-rw-310}
  \left| \partial_i \, \wtilde{p}^{n-w} (\t) \right| = (n-w) \left|
  \wtilde{p}^{n-w-1} (\t) \, \partial_i \wtilde{p}(\t) \right| \le
  \con n |\t| \le \con n^{(1+\eps)/2},
\end{equation}
via (\ref{t-rw-250}). As for the second sub-case, without loss of
generality, we assume the worst-case estimate for $\partial_i^{j_u}
\wtilde{p}$ on $\mathcal{B}_n$, that is,
\begin{equation}
  \label{t-rw-320}
  \left| \partial_i^{j_u} \wtilde{p} (\t) \right| \,\le\, \left\{
  \begin{array}{lll}
    \con |\t| \le \con n^{-(1-\eps)/2}, && \mbox{if $j_u = 1$,} \\
    \con, && \mbox{if $j_u \ge 2$.} \\ 
  \end{array}
  \right. 
\end{equation}
This implies that increasing by one the order of the derivative in the
l.h.s.\ of (\ref{t-rw-320}) will worsen our most conservative estimate
at most by a factor $n^{(1 - \eps)/2}$. Considering (\ref{t-rw-310})
as well, we conclude that applying another derivative to $\partial_i^l
\wtilde{p}^n$ will change its estimate at most by a factor $n^{(1 +
\eps)/2}$, which is again consistent with our inductive step.
\qed

\begin{remark}
  The careful reader might worry that the unrigorous use of the symbol
  $\con$ for a generic constant may jeopardize the above proof. It
  does not, since all the constants that have been used do not depend
  on $r$ or $n$. In principle, they may depend on $k$ (though it is
  easy to see that they do not), or $i$, or $j_u$, but these integers
  only take on a finite number of values, so bounds can be found that
  do not depend on any of the variables.
\end{remark}

\subsection{Stage 2: Removing (E3)}

If, contrary to assumption (E3), $v = -\imath \nabla
\wtilde{p} (0) \ne 0$, cf.\ (\ref{p3-20}), we define $\delta^{(n)} \in
\Z^d$ to be the (not necessarily unique) lattice point for which
$\delta^{(n)} / n$ best approximates $v \in \R^d$. One clearly has
\begin{equation}
  \label{t-rw-340}
  \left| \frac{\delta^{(n)}_i}{n} - v_i \right| \le \frac1{2n},
\end{equation}
where the subscript $i$ denotes, as usual, the $i$-th component of a
$d$-dimensional vector. Now, for $\t \in (-\pi, \pi)^d$, set
\begin{equation}
  \label{t-rw-350}
  \wtilde{\pi}_n (\t) := \wtilde{p} (\t) \, \emi{ (\delta^{(n)} \cdot 
  \t) / n}.
\end{equation}
We want to interpret $\wtilde{\pi}_n$ as a generally discontinuous
\fn\ $\T^d \longrightarrow \C$. On the other hand, $\wtilde{\pi}_n^n$
is a smooth \fn\ of $\T^d$ and, by (\ref{t-rw-140}) and Lemma
\ref{lem-anorm} of the Appendix,
\begin{eqnarray}
  \left\| \wtilde{g^{(n)}} \right\|_{\ac} &=& \left\| \left( 1 - 
  \wtilde{ q^{(r_n)} } \right) \wtilde{p}^n \right\|_{\ac} = 
  \nonumber \\
  &=& \left\| \left( 1 - \wtilde{ q^{(r_n)} } \right) \wtilde{p}^n 
  \, \wtilde{\omega}_{-\delta^{(n)}} \right\|_{\ac} = \nonumber \\
  \label{t-rw-360}
  &=& \left\| \left( 1 - \wtilde{q^{(r_n)} } \right) \wtilde{\pi}_n^n 
  \right\|_{\ac} 
\end{eqnarray}
(having used notation (\ref{def-phitilde}) as well).  Comparing the
above with (\ref{t-rw-150}), in view of our goal (\ref{t-rw-goal4}),
we see that it is sufficient to repeat all the estimations of Sections
\ref{subs-prop}--\ref{subs-endst1}, replacing $\wtilde{p}$ with
$\wtilde{\pi}_n$. This is no problem, except for estimate
(\ref{t-rw-250})---which is also reflected in (\ref{t-rw-310}) and
(\ref{t-rw-320}). (Consider Remark \ref{rk-no-ea3} and the fact that
Lemmas \ref{lem-pq2} and \ref{lem-p4} cannot distinguish between
$\wtilde{p}$ and $\wtilde{\pi}_n$.)

In order to find an effective substitute for (\ref{t-rw-250}), we
write, for $\t \in \mathcal{B}_n$,
\begin{equation}
  \label{t-rw-380}
  \partial_i \wtilde{\pi}_n (\t) = \partial_i \wtilde{\pi}_n
  (0) + u_n(\t') \cdot \t,
\end{equation}
where $u_n(\t')$ is the $i$-th row of the Hessian of $\wtilde{\pi}_n$
evaluated at some $\t' \in \mathcal{B}_n$. This has a finite limit, for $n
\to \infty$, as one can easily verify by direct computation on
(\ref{t-rw-350}) (it is in fact, up to a minus sing, the $i$-th row of
the covariance matrix of $p$). Therefore
\begin{equation}
  \label{t-rw-390}
  |u_n(\t') \cdot \t| \le \con |\t| \le \con n^{-(1-\eps)/2}.
\end{equation}
On the other hand, by (\ref{t-rw-340}),
\begin{equation}
  \label{t-rw-400}
  | \partial_i \wtilde{\pi}_n (0) | = \left| \partial_i
  \wtilde{p} (0) - \imath \frac{\delta^{(n)}_i}{n} \right| = \left| v_i
  - \frac{\delta^{(n)}_i}{n} \right| \le \con n^{-1}.
\end{equation}
Thus, using (\ref{t-rw-390})-(\ref{t-rw-400}) in (\ref{t-rw-380}),
\begin{equation}
  \label{t-rw-410}
  | \partial_i \wtilde{\pi}_n (\t) | \le \con n^{-(1-\eps)/2}, 
\end{equation}
which is the same bound as (\ref{t-rw-250}). This proves that the
$H^\nu$-norm of $\wtilde{g^{(n)}}$ can be estimated as in Section
\ref{subs-endst1} even when $\wtilde{p}$ is replaced by
$\wtilde{\pi}_n$. That is, (\ref{t-rw-goal}) holds even when
(E3) does not, which completes Stage 2.

\subsection{Stage 3: Removing (E1) and (E2)}

It is easy to realize that the convergence rate in (\ref{t-rw-goal})
is not only independent of the choice of $Q \subseteq S_0$, it is also
independent of the fact that $Q$ is contained in $S_0$, as long as it
remains an element of the countable partition associated to
$\scb_{-m,0}$.  In fact, if we take $Q \subseteq S_\a$, with $\a \ne
0$, we can shift, via the natural action of $\Z^d$ onto $\ps$, both
$Q$ and $F$ by the quantity $-\a$. (\ref{t-rw-goal}) continues to hold
with the same convergence rate because all the properties of $F$ that
were used in Stages 1 and 2 are translation invariant.

In formula, there exists a positive vanishing sequence $\{
\vartheta_n^{(m)} \}_{n \in \N}$ such that, if $g = 1_Q$ and $Q$ is a
fundamental set of $\scb_{-m,0}$,
\begin{equation}
  \label{t-rw-420}
  \left| \mu((F \circ T^n)g) - \avg(F) \mu(g) \right| \le \|g\|_{L^1} 
  \, \vartheta_n^{(m)} .
\end{equation}
Since (\ref{t-rw-420}) depends continuously on $g \in L^1$, it is
immediate to extend it to $g = \sum_{j \in \N} a_j 1_{Q_j}$, with $a_j
> 0$, that is, to a generic positive \fn\ in $L^1 (\ps, \scb_{-m,0},
\mu)$. If $g$ is such that both the positive part $g^+$ and the
negative part $g^-$ are nonzero, we apply (\ref{t-rw-420}) twice to
$g^+$ and $g^-$. An easy estimate proves that the formula holds in
this case as well.

Therefore \textbf{(M5)} holds w.r.t.\ $\what{\go} := \rset{F \in
L^\infty (\ps, \sca_s, \mu)} {\exists \, \avg(F)}$ and $\what{\lo}_m
:= L^1 (\ps, \scb_{-m,0}, \mu)$, for all $m \in \N$.

Now, if $F \in \go_m$ and $g \in \lo_m$, the invariance of $\mu$ and
Lemma \ref{lem-inv} give that
\begin{eqnarray}
  \label{t-rw-430}
  && \mu((F \circ T^n)g) - \avg(F) \mu(g) = \nonumber \\
  &=& \mu((F \circ T^{n-m}) (g \circ T^{-m})) - \avg(F \circ T^m) 
  \mu(g \circ T^{-m}).
\end{eqnarray}
Since $g \circ T^{-m} \in \what{\lo}_{2m}$ and $F \circ T^{m} \in
\what{\go}$ (because $\scb_{0,2m} \subset \sca_s$), we apply the
previous result and see that \emph{(a)} holds with a convergence rate
$\vartheta_{n-2m}^{(2m)}$.  

\subsection{Stage 4: Proof of the remaining assertions}

Statement \emph{(a)} immediately implies \textbf{(M4)} relative to
$\go_m$ and $\lo_m$ (Proposition \ref{prop-hye}). One readily extends
it to $\go$ and $\lo$, thus proving \emph{(b)}, by means of the
following obvious lemma:

\begin{lemma}
  If $\go'$ is a dense subset of $\go$ in the $L^\infty$-norm and
  $\lo'$ is a dense subset of $\lo$ in the $L^1$-norm, then
  \emph{\textbf{(M4)}} for $\go'$ and $\lo'$ implies
  \emph{\textbf{(M4)}} for $\go$ and $\lo$.
\end{lemma}

As concerns \emph{(c)}, it is easy to verify that Proposition
\ref{prop-m5m2} applies to the classes of global \ob s $\go_m$ and
local \ob s $\lo_m$ (using the family of local \ob s $g_\a := G \,
1_{S_\a}$). Therefore \emph{(a)} implies \textbf{(M2)} relative to
$\go_m$. We extend it to $\go$ by means of another obvious result.

\begin{lemma}
  If $\go'$ is a dense subset of $\go$ in the $L^\infty$-norm, then
  \emph{\textbf{(M2)}} for $\go'$ implies \emph{\textbf{(M2)}} for
  $\go$.
\end{lemma}

Finally, let us consider \emph{(d)}. By the second part of Proposition
\ref{prop-hye}, it suffices to show that, if $F, G \in \go_m$ and $F$
is $\Z^d$-periodic, then $\avg((F \circ T^n)G)$ exists for $n$ large
enough. By the same arguments as in the proof of Lemma \ref{lem-inv},
when $V \nearrow \ps$,
\begin{equation}
  \label{t-rw-440}
  \mu_V((F \circ T^n)G) = \mu_V((F \circ T^{n-m}) (G \circ T^{-m})) +
  o(1).
\end{equation}
So we can reduce to proving the existence of the infinite-volume limit
of the above r.h.s., for all $n \ge 2m$. Since $G \circ T^{-m}$ is
measurable w.r.t.\ $\scb_{-2m,0} \subset \sca_u$, with a slight abuse
of notation we can define, for $\a \in Z^d$,
\begin{equation}
  \label{t-rw-450}
  b_\a := \int_{S_\a} G \circ T^{-m} \, d\mu = \int_0^1 G \circ
  T^{-m} (\a; y_2) \, dy_2.
\end{equation}
An analogous definition can be made for $F \circ T^{n-m}$, which is
measurable w.r.t.\ $\scb_{0,2m} \subset \sca_s$. In this case, notice
that $F \circ T^{n-m}$ is also $\Z^d$-periodic, so we can write
\begin{equation}
  \label{t-rw-460}
  a := \int_{S_\a} F \circ T^{n-m} \, d\mu = \int_0^1 F \circ
  T^{n-m} (y_1) \, dy_1.
\end{equation}
Clearly, then,
\begin{equation}
  \label{t-rw-470}
  \int_{S_\a} (F \circ T^{n-m}) \, (G \circ T^{-m}) \, d\mu = a \, b_\a
\end{equation}
and, for $V = \bigcup_{\a \in B_{\g,r}} S_\a$,
\begin{equation}
  \label{t-rw-480}
  \mu_V((F \circ T^{n-m}) (G \circ T^{-m})) = \frac{a} {(2r+1)^d}
  \sum_{\a \in B_{\g,r}} b_\a.
\end{equation}
Since $\avg(G)$ exists, by Lemma \ref{lem-inv}
\begin{equation}
  \label{t-rw-490}
  \avg(G \circ T^{-m}) = \lim_{r \to \infty} \, \frac1 {(2r+1)^d}
  \sum_{\a \in B_{\g,r}} b_\a
\end{equation}
exists and the limit is uniform in $\g$. Also, it is obvious that
$\avg(F \circ T^{n-m}) = a$. Hence, as $V \nearrow \ps$, the r.h.s.\
of (\ref{t-rw-480}) tends to $\avg(F \circ T^{n-m}) \avg(G \circ
T^{-m})$, which is what we wanted to prove.

This concludes the proof of Theorem \ref{thm-rw}.
\qed

\appendix

\section{Appendix}

We collect here a few technical results which would have been
distracting in the body of the paper. The most important of them is an
estimate on a certain Fourier norm that is pivotal in the proof of
Theorem \ref{thm-rw}. This is presented in Section \ref{subs-anorm}.

\subsection{Proof of Lemma \ref{lem-a1}}
\label{subs-pf-a1}

Since $T$ is an automorphism, it is no loss of generality to prove the
assertion for $t = -1$. Also, since \textbf{(A1)} is invariant for the
action of $\Z^d$ on $V$, we may assume that all $V \in \scv$ are of
the form $V_r := B_{0,r} \times [0,1)^2$. Thus, the infinite-volume
limit becomes the limit $r \to \infty$.

Let $r' = r'(r) := [r^{1/2}]$ ($[ \cdot ]$ is the integer part of a
positive number) and
\begin{equation}
  \varphi(r') := \mu(T S_0 \setminus V_{r'}) = \sum_{\b \not\in
  B_{0,r'}} p_\b.
\end{equation}
Clearly, $\varphi(r') \searrow 0$, as $r'$ and $r$ tend to
infinity. This and the translation invariance of $T$ imply that
\begin{equation}
  \mu(T V_r \setminus V_{r+r'}) \le \mu(V_r) \, \varphi(r'),
\end{equation}
whence 
\begin{equation}
  \label{lem-a1-30}
  \mu(T V_r \cup V_r) \le \mu(V_{r+r'}) + \mu(T V_r \setminus
  V_{r+r'}) = (2r+1)^d + o((2r+1)^d).
\end{equation}
With a dual argument, considering that $\mu(V_{r-r'} \setminus T V_r)
= \mu(T^{-1} V_{r-r'} \setminus V_r)$ and that $T^{-1}$ acts
essentially as $T$ (after a swapping of the coordinates $y_1$ and
$y_2$, the map $T^{-1}$ becomes of the same type as $T$), we obtain
\begin{equation}
  \label{lem-a1-40}
  \mu(T V_r \cap V_r) \ge \mu(V_{r-r'}) - \mu(V_{r-r'} \setminus 
  T V_r) = (2r+1)^d + o((2r+1)^d).
\end{equation}
Taking the difference of (\ref{lem-a1-30}) and (\ref{lem-a1-40})
yields \textbf{(A1)} with $t = -1$.
\qed

\subsection{Proof of Lemma \ref{lem-span}}
\label{subs-pf-span}

We must show that, if $j',j'' \in \Z_N$ with $j' \ne j''$, then
\begin{equation}
  \label{l-span-10}
  \mathrm{span}_\Z \{ \b^{(j)} - \b^{(j')} \}_{j \ne j'} = 
  \mathrm{span}_\Z \{ \b^{(j)} - \b^{(j'')} \}_{j \ne j''}.
\end{equation}
The generic element of the l.h.s.\ of (\ref{l-span-10}) is 
\begin{equation}
  \label{l-span-20}
  \g = \sum_{j \ne j'} n_j (\b^{(j)} - \b^{(j')}) = \sum_{j \ne j'}
  n_j \b^{(j)} - \left( \sum_{j \ne j'} n_j \right) \b^{(j')},
\end{equation}
where $\{ n_j \}_{j \ne j'}$ are free variables, i.e., are arbitrarily
chosen integers. Upon defining $n_{j'} := - \sum_{j \ne j'} n_j$,
which implies $n_{j''} = - \sum_{j \ne j''} n_j$, (\ref{l-span-20})
becomes
\begin{equation}
  \label{l-span-30}
  \g = \sum_{j \in \Z_N} n_j \b^{(j)} = \sum_{j \ne j''} n_j \b^{(j)}
  - \left( \sum_{j \ne j''} n_j \right) \b^{(j'')} = \sum_{j \ne j''} 
  n_j (\b^{(j)} - \b^{(j'')}),
\end{equation}
which is the generic element of the r.h.s.\ of (\ref{l-span-10}), if
we consider $\{ n_j \}_{j \ne j''}$ to be the free variables and
$n_{j''}$ to depend on them.  
\qed

\subsection{Absolutely convergent Fourier series}
\label{subs-anorm}

In this section we present a convenient estimate for the space $\ac$
of \fn s $\wtilde{a}: \T^d \longrightarrow \C$ with an absolutely
convergent Fourier series $\{ a_\b \} = a$ \cite[\S6]{k}. This \fn al
space is defined as the maximal domain of the norm
\begin{equation}
  \| \wtilde{a} \|_{\ac} := \| a \|_{\ell^1} := \sum_{\b \in \Z^d} |a_\b|.
\end{equation}
This norm has a couple of straightforward invariances. For $\g \in
\Z^d$, $\zeta \in \T^d$, and $\wtilde{a} \in \ac$, let
\begin{eqnarray}
  \label{def-phitilde}
  \wtilde{\omega}_\g (\t) &:=& \ei{\g \cdot \t}; \\
  \label{def-tau}
  (\tau_\zeta \wtilde{a}) (\t) &:=& \wtilde{a} (\t + \zeta).
\end{eqnarray}

\begin{lemma}
  \label{lem-anorm}
  Given $\wtilde{a} \in A$, for all $\g \in \Z^d$ and $\zeta \in \T^d$,
  \begin{displaymath}
    \| \wtilde{a} \wtilde{\omega}_\g \|_{\ac} = \| \tau_\zeta \wtilde{a} 
    \|_{\ac} = \| \wtilde{a} \|_{\ac}.
  \end{displaymath}
\end{lemma}

\proof Trivial verification upon computation of the Fourier series of
$\wtilde{a} \wtilde{\omega}_\g$ and $\tau_\zeta \wtilde{a}$.
\qed

The following estimate is a modification---mostly, a
simplification---of a 1984 result by Nowak \cite{no}. The proof, which
we give for completeness, is practically copied from that article.

\begin{lemma}
  \label{lem-e-anorm}
  Let $\bar{\nu} = [d/2] + 1$ be the smallest integer strictly bigger
  than $d/2$. There exists a constant $C_d > 0$ such that
  \begin{displaymath}
    \| \wtilde{a} \|_{\ac} \le C_d \| \wtilde{a} \|_{H^{\bar{\nu}}},
  \end{displaymath}
  where
  \begin{eqnarray*}
    \| \wtilde{a} \|_{H^{\bar{\nu}}} &:=& |a_0| + \sum_{i=1}^d \left( 
    \sum_{\b \in \Z^d} \left| \b_i^{\bar{\nu}} a_\b \right|^2
    \right)^{\hspace{-4pt} 1/2} = \\
    &=& \left| \int_{\T^d} \wtilde{a} (\t) \dpt \right| + 
    \sum_{i=1}^d \left( \int_{\T^d} \left| \frac{ \partial^{\bar{\nu}}
    \wtilde{a} } { \partial \t_i^{\bar{\nu}} } (\t) \right|^2 
    \dpt \right)^{\hspace{-4pt} 1/2} .
  \end{eqnarray*}
\end{lemma}

\proof Let $\sigma = (\sigma_1, \sigma_2, \ldots, \sigma_d)$ be a
permutation of $(1,2, \ldots, d)$. Let us define 
\begin{equation}
  \label{anorm-10}
  Z_\sigma := \lset{(\b_1, \b_2, \ldots, \b_d) \in \Z^d}
  {0 < |\b_{\sigma_1}| \le |\b_{\sigma_2}| \cdots \le |\b_{\sigma_d}|}.
\end{equation}
Clearly, 
\begin{equation}
  \label{anorm-20}
  \Z^d = \{ 0 \} \cup \bigcup_\sigma Z_\sigma, 
\end{equation}
although the union is not disjoint. Using the Cauchy-Schwartz inequality,
\begin{equation}
  \label{anorm-30}
  \left( \sum_{\b \in Z_\sigma} |a_\b| \right)^2 \le C_{\bar{\nu}} 
  \sum_{\b \in Z_\sigma} \b_{\sigma_d}^{2\bar{\nu}} \, |a_\b|^2 \le
  C_{\bar{\nu}} \sum_{\b \in \Z^d} \left| \b_{\sigma_d}^{\bar{\nu}}
  a_\b\right|^2, 
\end{equation}
where we have denoted
\begin{eqnarray}
  C_{\bar{\nu}} := \sum_{\b \in Z_\sigma} \b_{\sigma_d}^{-2\bar{\nu}}
  \hspace{-6pt} &=& \hspace{-4pt} \sum_{|\a_1| > 0} \ 
  \sum_{|\a_2| \ge |\a_1|} \cdots \sum_{|\a_{d-1}| \ge |\a_{d-2}|} \ 
  \sum_{|\a_d| \ge |\a_{d-1}|} \a_d^{-2\bar{\nu}} \le \nonumber \\
  \label{anorm-40}
  &\le& \hspace{-4pt} \con \sum_{|\a_1| > 0} \ \sum_{|\a_2| \ge |\a_1|} 
  \cdots \sum_{|\a_{d-1}| \ge |\a_{d-2}|} \a_{d-1}^{-2\bar{\nu} + 1}
  \le \nonumber \\
  &\le& \cdots \cdots \hspace{8pt} \le \ \con 
  \sum_{|\a_1| > 0} \: \a_1^{-2\bar{\nu} + d - 1} < \infty
\end{eqnarray}
(as in Section \ref{sec-pf-main}, $\con$ represents a generic
constant). 
In view of (\ref{anorm-20}), summing the square root of
(\ref{anorm-30}) over all the permutations $\sigma$, we obtain
\begin{equation}
  \label{anorm-50}
  \sum_{\b \ne 0} |a_\b| \:\le\: (d-1)! \sqrt{C_{\bar{\nu}}} \:
  \sum_{i=1}^d  \left( \sum_{\b \in \Z^d} \left| \b_i^{\bar{\nu}} 
  a_\b \right|^2 \right)^{\hspace{-4pt} 1/2}.
\end{equation}
whence the assertion of the lemma.
\qed

\end{document}